\documentclass[a4paper,reqno,11pt]{amsart}

\usepackage[a4paper, scale={0.7,0.74}, marginratio={1:1}, footskip=7mm, headsep=10mm]{geometry}

\usepackage{amsmath,amsfonts,amssymb,amsthm,amscd,mathrsfs,mathdots,bbm,color,bm, stmaryrd}
\usepackage{mathtools}
\usepackage{graphicx}

\usepackage{xcolor}

\definecolor{linkblue}{RGB}{0,70,140}
\definecolor{citegreen}{RGB}{0,110,70}
\definecolor{urlpurple}{RGB}{100,60,130}

\usepackage[
    colorlinks=true,
    linkcolor=linkblue,
    citecolor=citegreen,
    urlcolor=urlpurple
]{hyperref}

\usepackage[normalem]{ulem}

\usepackage[style=alphabetic, backend=bibtex8, citestyle=alphabetic, giveninits=true, maxnames=99, maxalphanames=4, sorting=nyt]{biblatex}
\DeclareFieldFormat[article, inbook, incollection, inproceedings, patent, thesis, unpublished]{title}{#1} 
\renewbibmacro{in:}{}
\appto{\bibsetup}{\sloppy}

\makeatletter
\renewcommand\subsection{\@startsection{subsection}{2}%
  \z@{.5\linespacing\@plus.7\linespacing}{-.5em}%
  {\scshape}}
\makeatother

\def\XX{\mathcal{X}}
\def\R{{\mathbb R}}

\def\UT{\scalebox{0.58}{$\mathsf{UT}$}}
\def\P{{\mathbb P}}

\newcommand{\diff}{\mathrm{d}}

\newcommand{\Exp}{\mathrm{Exp}}
\newcommand{\Gumbel}{\mathrm{Gumbel}}

\newcommand{\botUT}{\bot^{\UT}}
\newcommand{\topUT}{\top^{\UT}}

\makeatletter
\DeclarePairedDelimiter\abs{\lvert}{\rvert} 
\let\oldabs\abs
\def\abs{\@ifstar{\oldabs}{\oldabs*}}
\DeclarePairedDelimiterX{\norm}[1]{\lVert}{\rVert}{#1} 
\let\oldnorm\norm
\def\norm{\@ifstar{\oldnorm}{\oldnorm*}}
\DeclarePairedDelimiterX{\ceil}[1]{\lceil}{\rceil}{#1} 
\let\oldceil\ceil
\def\ceil{\@ifstar{\oldceil}{\oldceil*}}
\DeclarePairedDelimiterX{\floor}[1]{\lfloor}{\rfloor}{#1} 
\let\oldfloor\floor
\def\floor{\@ifstar{\oldfloor}{\oldfloor*}}
\makeatother

\numberwithin{equation}{section}

\theoremstyle{plain}
\newtheorem{thm}{Theorem}[section]
\newtheorem{lem}[thm]{Lemma}
\newtheorem{prop}[thm]{Proposition}
\newtheorem{cor}[thm]{Corollary}

\theoremstyle{definition}
\newtheorem{definition}[thm]{Definition}

\newtheorem{rem}[thm]{Remark}
\newtheorem{ex}[thm]{Example}

\usepackage{etoolbox}
\AtEndEnvironment{rem}{\hfill$\triangleleft$}
\AtEndEnvironment{ex}{\hfill$\triangleleft$}

\newenvironment{myenumerate}{%
\renewcommand{\theenumi}{(\roman{enumi})}%
\renewcommand{\labelenumi}{\theenumi}%
\begin{list}{\labelenumi}
	{%
	\setlength{\itemsep}{0.4em}%
	\setlength{\topsep}{0.5em}%
	\setlength\leftmargin{2.45em}%
	\setlength\labelwidth{2.05em}%
	\setlength{\labelsep}{0.4em}%
	\usecounter{enumi}%
	}%
	}%
{\end{list}
}
\renewenvironment{enumerate}{
\begin{myenumerate}}%
{\end{myenumerate}}

\newcommand{\be}{\begin{equation}}
\newcommand{\ee}{\end{equation}}

\title[Extreme value statistics for majorisation partial orders]{Extreme value statistics \\ for majorisation partial orders}

\author[E.~Bisi]{Elia Bisi}
\address{Elia Bisi -- Universit\`a degli Studi di Firenze, Dipartimento di Matematica e Informatica, Viale Morgagni 67/A, 50134 Firenze, Italy}
\email{elia.bisi@unifi.it}

\author[F.~D.~Cunden]{Fabio Deelan Cunden}
\address{Fabio Deelan Cunden -- Dipartimento di Matematica, Universit\`a degli Studi di Bari, I-70125 Bari, Italy, and
INFN, Sezione di Bari, I-70126 Bari, Italy}
\email{fabio.cunden@uniba.it}

\author[G.~Gramegna]{Giovanni Gramegna}
\address{Giovanni Gramegna -- Dipartimento di Fisica, Universit\`a degli Studi di Bari, I-70165 Bari, Italy, and
INFN, Sezione di Bari, I-70126 Bari, Italy}
\email{giovanni.gramegna@uniba.it}

\author[M.~Ligab\`o]{Marilena Ligab\`o}
\address{Marilena Ligab\`o -- Dipartimento di Matematica, Universit\`a degli Studi di Bari, I-70125 Bari, Italy}
\email{marilena.ligabo@uniba.it}


\begin{document}

\begin{abstract}
We study extreme-value phenomena for random samples in partially ordered sets through their meets and joins, the order-theoretic analogues of minima and maxima. In this paper we consider the probability simplex equipped with either UT-majorisation or majorisation, with samples drawn uniformly from the corresponding simplex.
We analyse these extrema both at fixed dimension and as the dimension grows, identifying almost-sure collapse to the bottom or top elements and critical scalings with nondegenerate limiting laws.
\end{abstract}

\maketitle

\tableofcontents

\section{Introduction}

The study of extreme values of real random variables and their limiting laws is a cornerstone of probability theory.
Classical questions concern the asymptotic behaviour, as $m$ tends to infinity, of $\min\{X_1,\ldots,X_m\}$ and $\max\{X_1,\ldots,X_m\}$, where $X_1,X_2,\dots$ are i.i.d.\ real random variables and the extrema are taken with respect to the usual total order $\leq$ on $\R$.
Starting from the problem of determining when such extrema converge after suitable normalisation, a vast theory has developed~\cite{Gumbel58,Resnick87,Majumdar2024}.
Applications of extreme value theory range widely, from climate science to finance; see, for example,~\cite[Chap.~1]{Beirlant2004}.

The real line is, however, a special setting: it is one-dimensional and totally ordered.
The extension of extreme value theory beyond the real line has a long history.
Early work by Gumbel~\cite{Gumbel1962} already considered multivariate extremal distributions based on componentwise minima and maxima.
This approach has developed into a rich multivariate extreme value theory; see, for instance,~\cite[Chaps.~8--9]{Beirlant2004} and the references therein.
Classical multivariate extreme value theory typically retains the coordinate or geometric structure of Euclidean space, studying extremal behaviour in componentwise, directional or radial terms.

More generally, however, there is no canonical notion of an extreme observation in a multidimensional space, since such spaces do not generally come equipped with a  total order.
In many settings, including random graphs, matrices, integer partitions and probability vectors, the natural objects are both high-dimensional and equipped with a partial rather than a total order.
Most pairs of elements may then be incomparable, so the usual notions of sample maximum and sample minimum no longer directly apply and must be suitably generalised.

This article studies extreme-value phenomena for random samples in partially ordered sets through the asymptotic behaviour of their meets and joins, the lattice-theoretic analogues of minima and maxima.
Our perspective places the order structure at the centre of the analysis, defining extremes intrinsically through a prescribed partial order.
This framework is particularly suited to settings in which incomparability is a structural feature rather than a technical complication.
In particular, we consider UT-majorisation on the $(n-1)$-dimensional probability simplex and majorisation on the corresponding ordered simplex.
In both cases, the samples consist of independent random probability vectors drawn uniformly from the relevant simplex.

Our principal motivation comes from resource theories~\cite{Alberti82,winter2016operational,Chitambar19,deoliveirajunior2022}, where convertibility under allowed operations induces a partial order on states rather than a total ranking by resource content.
Majorisation and UT-majorisation arise naturally in several classical and quantum resource theories~\cite{Alberti82,Nielsen1999,Nielsen2001,winter2016operational,Chitambar19,vomEnde}.
In resource theories with a lattice structure, the meet of two states can be interpreted as an \emph{optimal common resource} state, while the join is sometimes referred to as an \emph{optimal common product} state.
See the recent papers~\cite{Deside2024,Massri2024} for the operational interpretation of the meet and the join in the resource theory of bipartite entanglement.

\section{Preliminaries}
\label{sec:preliminaries}

As a general introduction on partial orders and majorisation relations the reader can consult~\cite[Ch.~3]{Stanley2011} and~\cite[Part I]{MOA79}, respectively.

\subsection{Order-theoretic preliminaries}
A \textbf{partially ordered set (poset)} $\left(\XX,\sqsubseteq\right)$ is a set $\XX$ endowed with a \textbf{partial order}, i.e.\ a relation on $\XX$ that is reflexive, antisymmetric and transitive.

The \textbf{meet} or \textbf{infimum} or \textbf{greatest lower bound} of  $x,y \in \XX$ is, if it exists, the unique element $z\in\XX$ such that:
\begin{enumerate}
    \item $z\sqsubseteq x$ and $z\sqsubseteq y$;
    \item if $w\in\XX$ satisfies $w\sqsubseteq x$ and $w\sqsubseteq y$, then $w \sqsubseteq z$.
\end{enumerate}
Whenever it exists, the \textbf{join} or \textbf{supremum} or \textbf{least upper bound} of $x,y\in \XX$ is defined as the meet of $x,y$ for the \textbf{inverse order} $\sqsupseteq$ (i.e.\ the partial order defined by setting $x\sqsupseteq y$ if and only if $y\sqsubseteq x$).

In the case of a totally ordered set (for example $\left(\XX,\sqsubseteq\right)=(\R,\leq)$, the real numbers with the usual order), the meet and the join of $x,y$ are the \textbf{minimum} and \textbf{maximum} between $x$ and $y$, respectively; one of them is $x$ and the other is $y$.

A poset $\left(\XX,\sqsubseteq\right)$ is a \textbf{lattice} if every pair of elements $x,y\in\XX$ has a meet, denoted $x\wedge y$, and a join, denoted $x\vee y$. 
In a lattice, $\wedge$ and $\vee$ define two commutative ($x\wedge y=y\wedge x$), associative ($(x\wedge y)\wedge z=x\wedge (y\wedge z)$) and idempotent ($x\wedge x=x$) binary operations on $\XX$.
Therefore, we may unambiguously 
write
\[
\bigwedge_{i=1}^m x_i:=x_1\wedge x_2\wedge\cdots\wedge x_m \qquad\text{and}\qquad
\bigvee_{i=1}^m x_i:=x_1\vee x_2\vee\cdots\vee x_m
\]
for the meet and join, respectively, of $m$ elements $x_1,\ldots,x_m\in\XX$, without caring about the order in which the operations are performed. 
A lattice $(\XX,\sqsubseteq)$  is \textbf{bounded} if there exist two elements $\bot,\top\in\XX$ such that 
$\bot \sqsubseteq x\sqsubseteq \top$, for all $x\in\XX$. If they exist, such elements are necessarily unique. We call $\bot$ the \textbf{bottom} (element) and $\top$ the \textbf{top} (element) of the lattice.

\subsection{Majorisation orders on the simplex}
Let $n\geq 2$.
Throughout this paper, we will work with two partial orders on probability simplices.
The first is \textbf{UT-majorisation} (short for `upper-triangular majorisation'), a partial order on the unit simplex
\begin{equation}
\Delta_{n-1} := \{x=(x_1,\dots,x_n)\in \R^n_{\geq0}\colon x_1+\dots+x_n=1\}. 
\end{equation}
The second is \textbf{majorisation}, regarded as a partial order on the subset of the unit simplex given by the ordered Weyl chamber
\begin{equation}
\label{eq:WeylChamber}
\Delta_{n-1}^{\downarrow} := \{(x_1,\dots,x_n)\in \Delta_{n-1}\colon  x_1\geq x_2\geq\cdots\geq x_n\}
= \{x^{\downarrow}\colon x \in \Delta_{n-1}\} 
\end{equation}
(here and throughout, $x^\downarrow$ denotes the rearrangement of a vector $x$ in non-increasing order).
Equivalently, majorisation can be defined on the space of probability $n$-vectors modulo coordinate permutations.

We use the standard notation $\llbracket n\rrbracket := \{1,\ldots,n\}$. 
Let
\begin{equation}
\phi\colon \R^{n}\to\R^{n},\quad x=(x_1,\ldots,x_n)\longmapsto \phi(x)=\left(\phi_1(x),\ldots,\phi_n(x)\right), 
\end{equation}
where $\phi_k(x):=x_1+\cdots+x_k$ is the \textbf{$k$-th partial sum} of $x$, for all $k\in\llbracket n\rrbracket$. The linear map $\phi$ is a volume-preserving bijection with inverse 
\begin{equation}
\label{eq:phi-inverse}
\phi^{-1}(u)=
\left(u_1,u_2-u_1,\ldots,u_n-u_{n-1}\right).
\end{equation}
For later use, we also set $\phi_0(x):=0$ and define the \textbf{augmented partial-sum vector}
\begin{equation}
\label{eq:phi-hat}
\widehat{\phi}(x)
:=
\left(\phi_0(x),\phi_1(x),\ldots,\phi_n(x)\right)
=
\left(0,\phi_1(x),\ldots,\phi_n(x)\right).
\end{equation}

\begin{definition}[UT-majorisation]
\label{def:UT-majorisation}
    If $x,y\in \Delta_{n-1}$, we write $x \preceq  y$, and say that $x$ is \textbf{UT-majorised} by $y$, if $\phi_k(x)\leq \phi_k(y)$ for all $k\in \llbracket n\rrbracket$.
\end{definition}
\begin{definition}[Majorisation]
\label{def:Majorisation}
    If $x,y\in \Delta_{n-1}^{\downarrow}$, we write $x\preceq y$, and say that $x$ is \textbf{majorised} by $y$, if $\phi_k(x)\leq \phi_k(y)$ for all $k\in \llbracket n\rrbracket$.
\end{definition}
Note that $\phi_n(x)=\phi_n(y)=1$ for all $x,y\in\Delta_{n-1}$, hence the inequality for $k=n$ in the above definitions is redundant.

\begin{rem}
In most of the literature, majorisation is defined as a relation on the
entire unit simplex by comparing the decreasing rearrangements
$x^\downarrow$ and $y^\downarrow$ of $x,y\in\Delta_{n-1}$. Consequently,
vectors that differ only by a permutation majorise one another, so the
resulting relation is a preorder rather than a partial order. 
Every permutation class contains a unique element of
$\Delta_{n-1}^{\downarrow}$. Under Definition~\ref{def:Majorisation},
we work with these canonical representatives, and majorisation becomes
a \emph{bona fide} partial order. Moreover, it coincides with the
restriction of UT-majorisation to the Weyl chamber
\eqref{eq:WeylChamber}. We may therefore use the same symbol
$\preceq$ for both relations.
\end{rem}

It is easy to see that $(\Delta_{n-1},\preceq)$ and $(\Delta_{n-1}^{\downarrow},\preceq)$ are not totally ordered if $n\geq3$: there are probability vectors that are incomparable with respect to both UT-majorisation and majorisation. 

The following structural fact is the order-theoretic basis for the extreme-value questions studied in this paper.
\begin{thm}
\label{thm:lattices}
    The posets $(\Delta_{n-1},\preceq)$ and $(\Delta_{n-1}^{\downarrow},\preceq)$ are bounded lattices, with bottom and top elements defined by 
\begin{align}
\label{eq:UT_top_bottom}
    \botUT_n:=(0,\ldots,0,1)&\preceq x\preceq (1,0,\ldots,0)=:\topUT_n &&\text{for all $x\in\Delta_{n-1}$},\\
    \label{eq:maj_top_bottom}
        \bot_n:=\left(\frac{1}{n},\ldots,\frac{1}{n}\right)&\preceq x\preceq(1,0,\ldots,0)=:\top_n &&\text{for all $x\in\Delta_{n-1}^{\downarrow}$}.
\end{align}
\end{thm} 
Both lattices admit explicit descriptions of their meet and join. We first introduce some notation.

For $a=(a_1,\ldots,a_n)$ and $b=(b_1,\ldots,b_n)$ in $\R^n$, write $a\leq b$ when $a_k\leq b_k$ for every $k\in\llbracket n\rrbracket$, and let 
\begin{align*}
\min(a,b)&:=(\min(a_1,b_1),\ldots,\min(a_n,b_n)), \\ \max(a,b)&:=(\max(a_1,b_1),\ldots,\max(a_n,b_n))
\end{align*}
their componentwise minimum and maximum.

Following~\cite[p.~11]{Mitrinovic70}, we say that a real tuple $a=(a_0,\ldots,a_{n})$ is \textbf{concave} if 
\[
2a_k\geq a_{k-1}+a_{k+1},\quad \text{for all $k\in \llbracket n-1 \rrbracket$}.    
\]

For $a\in\R^n$, we denote by $a^*$ its \textbf{concave envelope}, namely the smallest concave tuple lying above $a$:
\begin{enumerate}
    \item $a^*$ is concave and $a\leq a^*$;
    \item $a^*\leq b$ whenever $b$ is a concave sequence such that $a\leq b$.
\end{enumerate}
The UT-majorisation meet and join are simple to describe.
\begin{prop}
\label{prop:meet-UT}
The meet $x\wedge y$ and join $x\vee y$ of $x,y\in\Delta_{n-1}$ with respect to the UT-majorisation relation are determined by
\begin{equation}
\label{eq:meet-join-UT}
    \phi(x\wedge y)=\min\left(\phi(x),\phi(y)\right),
        \qquad \phi(x\vee y)=\max\left(\phi(x),\phi(y)\right). 
\end{equation}
\end{prop}

For majorisation, the construction of the join requires an additional concavification step.
The augmented partial sums of every vector in the ordered simplex form a concave sequence.
The componentwise minimum of two concave sequences is again concave, so it directly determines the meet.
By contrast, their componentwise maximum need not be concave, and the partial-sum sequence of the join is obtained by replacing it with its concave envelope.

\begin{prop}
\label{prop:meet-majorisation}
The meet $x\wedge y$ and join $x\vee y$ of $x,y\in\Delta_{n-1}^{\downarrow}$ with respect to the majorisation relation are determined by
\begin{equation}
\label{eq:meet-join-maj}
    \phi(x\wedge y)=\min\left(\phi(x),\phi(y)\right), \qquad
        \widehat{\phi}(x\vee y)=\max\left(\widehat{\phi}(x),\widehat{\phi}(y)\right)^*. 
\end{equation}
\end{prop}
In equivalent notation, formulae~\eqref{eq:meet-join-UT} for UT-majorisation appear in~\cite{ParkerRam99}.
Formulae~\eqref{eq:meet-join-maj} for majorisation can be found in~\cite{Cicalese2002}, where the join is constructed algorithmically.
The formulation of the join through the concave envelope makes the underlying structure particularly transparent; to our knowledge, this explicit formulation has not been emphasised in the literature.

We prove Propositions~\ref{prop:meet-UT}--\ref{prop:meet-majorisation}, together with Theorem~\ref{thm:lattices}, in Appendix~\ref{sec:geometry}.

\section{Overview of the results}
\label{sec:overview}
What began as an exploration of extreme-value theory on majorisation lattices revealed a richer picture than we had anticipated.
The behaviour of meets and joins differs markedly between the two orders and between fixed and growing dimension.
We therefore begin with an informal overview of the main phenomena, leaving the precise statements to the next section.

Throughout the paper, unless explicitly stated otherwise, we let $\left(X_{n,i}\right)_{n\geq 2, i\geq1}$ be a collection of independent random vectors such that $X_{n,i}$ is uniformly distributed on $\Delta_{n-1}$.
Similarly, we let $\left(Y_{n,i}\right)_{n\geq 2, i\geq1}$ be a collection of independent random vectors, with $Y_{n,i}$ uniformly distributed on $\Delta_{n-1}^{\downarrow}$.
These samples may be realised on a common probability space in terms of i.i.d.\ uniform random variables; see Propositions~\ref{prop:UT-uniform} and~\ref{prop:maj-uniform}.

For all $n\geq 2$ and $m\geq1$,  denote the UT-majorisation meets and joins by
\begin{equation}
    \label{eq:UTmaj-meet-join}
\mathsf{\Lambda}_{n,m}^{\UT}:=\bigwedge_{i=1}^m X_{n,i} \qquad\text{and}\qquad \mathsf{V}_{n,m}^{\UT}:=\bigvee_{i=1}^m X_{n,i},
\end{equation}
and, similarly, the majorisation meets and joins by
 \begin{equation}
    \label{eq:maj-meet-join}
\mathsf{\Lambda}_{n,m}:=\bigwedge_{i=1}^m Y_{n,i} \qquad\text{and}\qquad \mathsf{V}_{n,m}:=\bigvee_{i=1}^m Y_{n,i}.
\end{equation}

We study the asymptotic behaviour of these random meets and joins in two regimes: fixed dimension and growing dimension.

\noindent
\textbf{Collapse at fixed dimension.}
For fixed $n$, both the meet and the join converge almost surely to the bottom and top elements of the respective order as $m\to\infty$. The corresponding results are Theorem~\ref{thm:UT-conv-as} for UT-majorisation and Theorem~\ref{thm:maj-conv-as} for majorisation.

\noindent
\textbf{Fluctuations at fixed dimension.}
The fluctuation scales reflect the geometry of the two orders.

For UT-majorisation, the $k$-th partial sum of the meet fluctuates on the scale $m^{-1/k}$, $k=1,\ldots,n-1$.
After the corresponding rescaling, the limiting coordinates are independent random variables whose $k$-th powers are exponentially distributed, with explicit rates.
By reflection, an analogous statement holds for the join near the top element; see Theorem~\ref{thm:UT-conv-distr}.

For majorisation, the meet behaves differently.
The concavity constraint imposed by the ordered simplex forces all partial-sum deviations from the bottom element onto the common scale $m^{-1/(n-1)}$.
The limiting vector is supported on the cone of nonnegative concave sequences and admits an explicit description in terms of volumes of truncations of this cone.
Its coordinates are generally dependent, although their one-dimensional marginals are again of exponential-power type; moreover, the limiting law is palindromic in distribution.
See Theorem~\ref{thm:maj-conv-distr} and Corollary~\ref{cor:maj-limit-marginals}.

The majorisation join has a different fluctuation regime.
For $k=1,\ldots,n-1$, the quantity $1-\phi_k(\mathsf{V}_{n,m})$ fluctuates on the scale $m^{-1/(n-k)}$.
The rescaled random variables converge jointly to independent random variables whose $(n-k)$-th powers are exponentially distributed, again with explicit  rates.
Thus, the meet has a common fluctuation scale but dependent limiting coordinates, whereas the join has coordinate-dependent fluctuation scales but independent limiting coordinates.
Moreover, the concavification involved in constructing the join (see Proposition~\ref{prop:meet-majorisation}) is asymptotically inactive; see Theorem~\ref{thm:maj-join-conv-distr}.

\noindent
\textbf{Growing dimension.}
We finally let $n\to\infty$ together with $m=m_n$.
Here the question is how rapidly the sample size must grow to overcome the concentration induced by the increasing dimension.

For UT-majorisation, the sup-norm distances of the meet from the bottom and of the join from the top are both almost surely asymptotic to $m_n^{-1/(n-1)}$.
If $\log m_n/n\to c\in[0,\infty]$, these distances converge almost surely to $e^{-c}$, with the convention $e^{-\infty}:=0$.
Consequently, convergence to the corresponding extremal elements occurs precisely when
\[
\frac{\log m_n}{n}\longrightarrow\infty.
\]
When $m_n\to\infty$, both distances have Gumbel fluctuations at the next order; see Theorem~\ref{thm:UT-m-n-scaling} and Corollary~\ref{cor:UT-m-n-collapse}.

For majorisation, the sup-norm distance of the join from the top is almost surely asymptotic to $(nm_n)^{-1/(n-1)}$ and follows the same $e^{-c}$ phase diagram.
Its next-order Gumbel limit holds without requiring $m_n\to\infty$, even when $m_n=1$.
The same superexponential condition characterises convergence of the join to the top and convergence of the meet to the bottom at the finer $o(1/n)$ scale in sup-norm.
See Theorems~\ref{thm:maj-meet-m-n-collapse} and~\ref{thm:maj-join-m-n-scaling}, and Corollary~\ref{cor:maj-join-m-n-collapse}.

\subsection{Organisation of the paper}
Section~\ref{sec:results} states the precise results for UT-majorisation and majorisation.
Section~\ref{sec:distr_identities} develops the distributional identities underlying the analysis, including the representation of uniform points on the simplex and ordered simplex in terms of spacings of independent uniform random variables, together with related Markov properties.
Section~\ref{sec:proofs} contains the proofs of the main results.
Section~\ref{sec:future} discusses open problems and possible extensions.
Finally, in Appendix~\ref{sec:geometry} we include the proofs of the explicit formulae for the meet and the join in the UT-majorisation and majorisation orders stated in Section~\ref{sec:preliminaries} and present their geometric interpretations in terms of Lorenz curves.

\section{Detailed results}
\label{sec:results}

We now state the main results of the paper.
The proofs are given in Sec.~\ref{sec:proofs_UTmaj} for UT-majorisation and in Sec.~\ref{sec:proofs_maj} for majorisation.

\subsection{Notation}

We write $\xi\sim\rho$ to indicate that the random variable $\xi$ has distribution $\rho$.
The symbol $\stackrel{\mathrm{d}}{=}$ denotes equality in distribution.
Almost sure convergence and convergence in distribution are denoted, respectively, by $\xrightarrow{\operatorname{a.s.}}$ and $\xrightarrow{\operatorname{d}}$.

Recall that an exponential random variable $\xi$ with rate $\lambda>0$ has distribution function
\begin{equation}
    \P(\xi\leq x)=
    \begin{cases}
        1-\exp(-\lambda x), & \text{if } x\geq 0,\\
        0, & \text{if } x<0.
    \end{cases}
\end{equation}
We write $\xi\sim\Exp(\lambda)$ in this case.

Recall also that a standard Gumbel random variable $G$ has distribution function
\begin{equation}
    \P(G\leq x)=\exp(-e^{-x}),
    \qquad x\in\R.
\end{equation}
We write $G\sim \Gumbel$ in this case.

\subsection{Fixed-dimensional results}

\begin{figure}
\centering
\includegraphics[width=\textwidth]{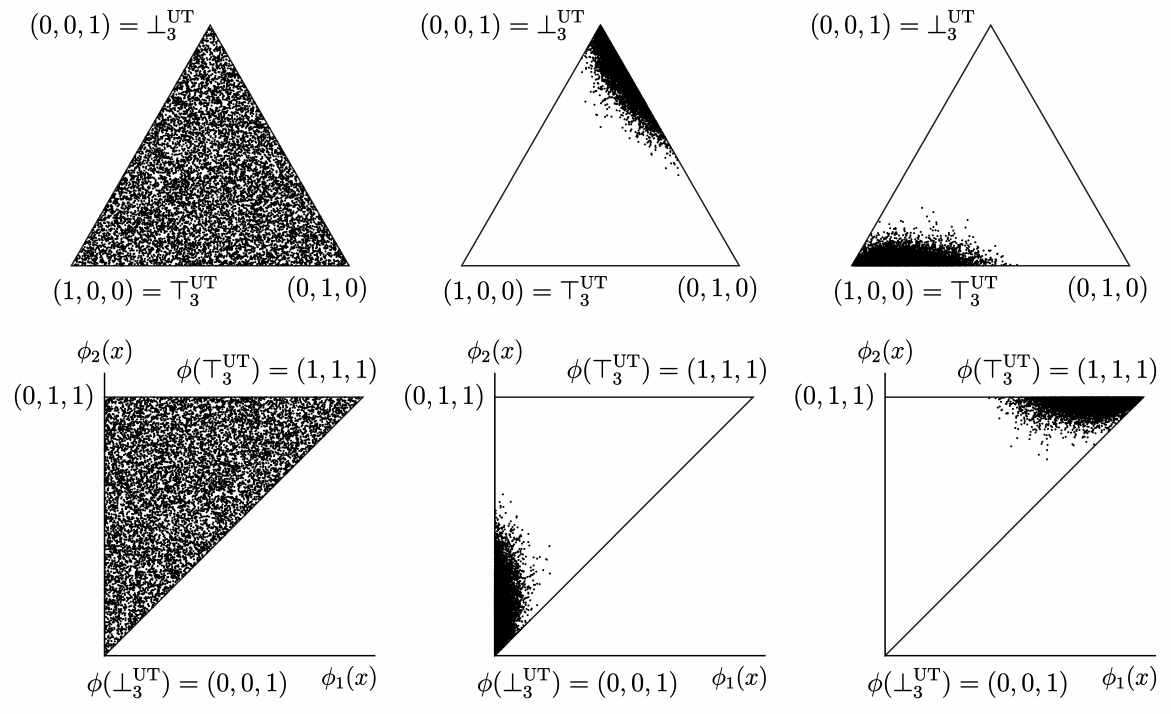}
\caption{Top row: numerical sample of $10^4$ random points on the simplex $\Delta_2$ ($n=3$); from left to right, the panels show the uniform distribution, the UT-majorisation meet of $m$ i.i.d.\ uniform points, and the UT-majorisation join of $m$ i.i.d.\ uniform points, with $m=20$.  
Bottom row: the same samples represented in the space of the partial sums.}
\label{fig:inf_simplex_UTmaj}
\end{figure}

\begin{figure}
\centering
\includegraphics[width=\textwidth]{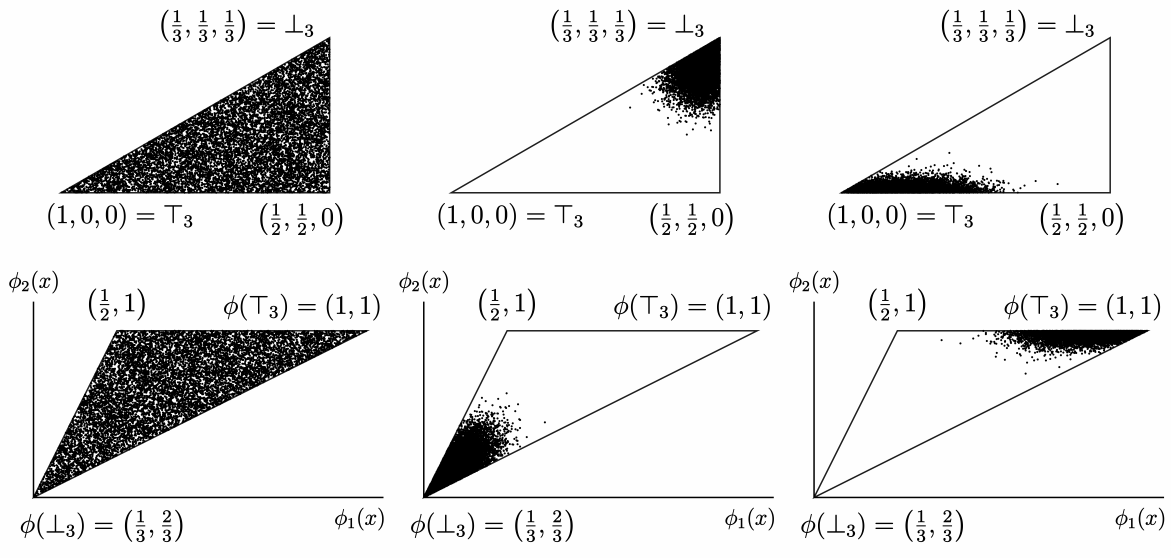}
\caption{
Top row: numerical sample of $10^4$ random points on the ordered simplex $\Delta_2^{\downarrow}$ ($n=3$); from left to right, the panels show the uniform distribution, the majorisation meet of $m$ i.i.d.\ uniform points, and the majorisation join of $m$ i.i.d.\ uniform points, with $m=20$.  
Bottom row: the same samples in the space of partial sums.}
\label{fig:inf_simplex_maj}
\end{figure}

Our first result concerns the almost sure convergence, for fixed space dimension $n$, of meet and join as the sample size $m$ grows large. 
The results are illustrated in Fig.~\ref{fig:inf_simplex_UTmaj} and~\ref{fig:inf_simplex_maj}.

\begin{thm}[Fixed-dimensional convergence of the UT-majorisation meet and join]
\label{thm:UT-conv-as}
Fix $n\geq2$ and let $\mathsf{\Lambda}_{n,m}^{\UT}$ and $\mathsf{V}_{n,m}^{\UT}$ be as in~\eqref{eq:UTmaj-meet-join}.
Then, as $m\to\infty$, $\mathsf{\Lambda}_{n,m}^{\UT}\xrightarrow{\operatorname{a.s.}}\botUT_n$ and $\mathsf{V}_{n,m}^{\UT}\xrightarrow{\operatorname{a.s.}}\topUT_n$.
\end{thm}

\begin{thm}[Fixed-dimensional convergence of the majorisation meet and join]
\label{thm:maj-conv-as}
Let $n\geq2$ and let $\mathsf{\Lambda}_{n,m}$ and $\mathsf{V}_{n,m}$ be as in~\eqref{eq:maj-meet-join}.
Then, as $m\to\infty$, $\mathsf{\Lambda}_{n,m}\xrightarrow{\operatorname{a.s.}}\bot_n$ and $\mathsf{V}_{n,m} \xrightarrow{\operatorname{a.s.}}\top_n$. 
\end{thm}

\begin{figure}
\centering
\includegraphics[width=0.7\textwidth]{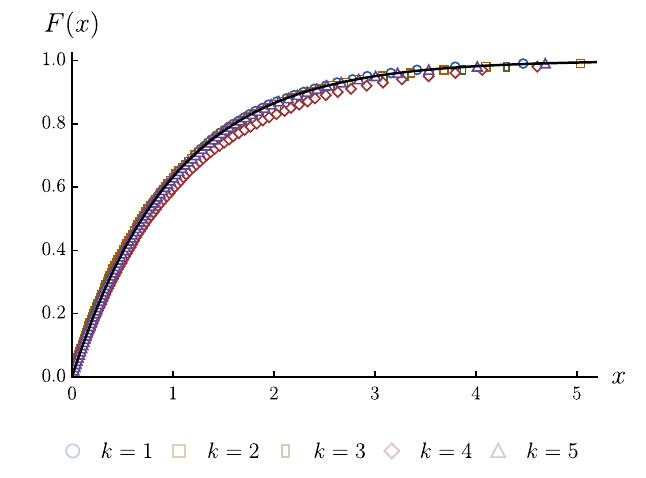}
\caption{Numerical illustration of the UT-majorisation scaling limit.
Empirical distribution function, computed from $10^4$ realisations, of $\lambda_{n,k}m\phi_k(\mathsf{\Lambda}_{n,m}^{\UT})^k$, with $\lambda_{n,k}$ given in~\eqref{eq:UT_conv_distr}.
The simulation is performed for $n=6$ and $m=1000$.
The limiting exponential distribution function with rate $1$ is shown as a solid black line.}
\label{fig:meet_UT_num}
\end{figure}

For a finer result, we study the scaling under which we have a non-trivial convergence in distribution as $m\to\infty$.
The theorem is illustrated in  Fig.~\ref{fig:meet_UT_num}.
\begin{thm}[Fixed-dimensional scaling limit for UT-majorisation]
\label{thm:UT-conv-distr}
Fix $n\geq2$ and let $\mathsf{\Lambda}_{n,m}^{\UT}$ and $\mathsf{V}_{n,m}^{\UT}$ be as in~\eqref{eq:UTmaj-meet-join}.
Then, as $m\to\infty$, we have
\begin{align}
\left(m^{\frac{1}{k}} \phi_k\left(\mathsf{\Lambda}_{n,m}^{\UT}\right) \right)_{k=1}^{n-1}
&\xlongrightarrow{\operatorname{d}} \left(\xi_1,\ldots,\xi_{n-1}\right) \\
\intertext{and}
\left(m^{\frac{1}{n-k}}(1-\phi_k\left(\mathsf{V}_{n,m}^{\UT}\right))\right)_{k=1}^{n-1} &\xlongrightarrow{\operatorname{d}}\left(\xi_{n-1},\ldots,\xi_1\right),
\end{align}
where the random variables $(\xi_k)_{k=1}^{n-1}$ are independent and 
\begin{equation}
\label{eq:UT_conv_distr}
   \xi_k^{k}\sim \Exp(\lambda_{n,k}),\qquad  \lambda_{n,k}=\binom{n-1}{k}.
\end{equation}
\end{thm}

\begin{rem}\label{rem:ScalingCoordinatesUT}
The convergence in Theorem~\ref{thm:UT-conv-distr} can equivalently be
expressed in the standard coordinates of the simplex.
Using~\eqref{eq:phi-inverse} and the fact that the successive partial sums fluctuate on different scales, we obtain, for $2\leq k\leq n-1$,
\begin{equation*}
    m^{\frac{1}{k}}(\mathsf{\Lambda}_{n,m}^{\UT})_k
=
m^{\frac{1}{k}}
\bigl(
\phi_k(\mathsf{\Lambda}_{n,m}^{\UT})
-\phi_{k-1}(\mathsf{\Lambda}_{n,m}^{\UT})
\bigr) \xlongrightarrow{\operatorname{d}} \xi_k.
\end{equation*}
Indeed, $\phi_{k-1}(\mathsf{\Lambda}_{n,m}^{\UT})$ is of order
$m^{-1/(k-1)}$ and, hence, $m^{1/k}\phi_{k-1}(\mathsf{\Lambda}_{n,m}^{\UT}) \xlongrightarrow{\operatorname{d}}0$.
Moreover, 
\begin{equation*}
    m^{\frac{1}{n-1}}\bigl(1-(\mathsf{\Lambda}_{n,m}^{\UT})_n\bigr)
= m^{\frac{1}{n-1}} \phi_{n-1}(\mathsf{\Lambda}_{n,m}^{\UT}) \xlongrightarrow{\operatorname{d}} \xi_{n-1}.
\end{equation*}
Therefore, we obtain
\begin{equation*}
    m^{\frac{1}{k}}(\mathsf{\Lambda}_{n,m}^{\UT})_k \xlongrightarrow{\operatorname{d}} \xi_k, \quad k\in\llbracket n-1\rrbracket, \qquad\text{and}\qquad
m^{\frac{1}{n-1}}\bigl(1-(\mathsf{\Lambda}_{n,m}^{\UT})_n\bigr) \xlongrightarrow{\operatorname{d}} \xi_{n-1}
\end{equation*}
jointly.
For the join, similarly, we have
\begin{equation*}
    m^{\frac{1}{k}}(\mathsf{V}_{n,m}^{\UT})_{n-k+1} \xlongrightarrow{\operatorname{d}} \xi_k, \quad k\in\llbracket n-1\rrbracket, \qquad\text{and}\qquad
m^{\frac{1}{n-1}}\bigl(1-(\mathsf{V}_{n,m}^{\UT})_1\bigr) \xlongrightarrow{\operatorname{d}} \xi_{n-1}
\end{equation*}
jointly.
\end{rem}

For UT-majorisation, the reflection symmetry 
\[
\left(\phi_k\left(\mathsf{\Lambda}_{n,m}^{\UT}\right)\right)_{k=1}^{n-1}\stackrel{\mathrm{d}}{=}\left(1-\phi_{n-k}\left(\mathsf{V}_{n,m}^{\UT}\right)\right)_{k=1}^{n-1}
\]
holds.
Equivalently, in standard coordinates,
\[    \mathsf{\Lambda}_{n,m}^{\UT}\stackrel{\mathrm{d}}{=} R(\mathsf{V}_{n,m}^{\UT}),
\]
where
\begin{equation}\label{eq:Reflection}
    R(x_1,\dots,x_n):=(x_n,...,x_1).
\end{equation}
Consequently, the scaling limits of the UT-majorisation meet and join can be deduced from one another.
Coordinate reversal does not preserve the ordered simplex; for $n\geq3$, the majorisation meet and join instead exhibit different fixed-dimensional fluctuation regimes.

\begin{figure}
\centering
\includegraphics[width=\textwidth]{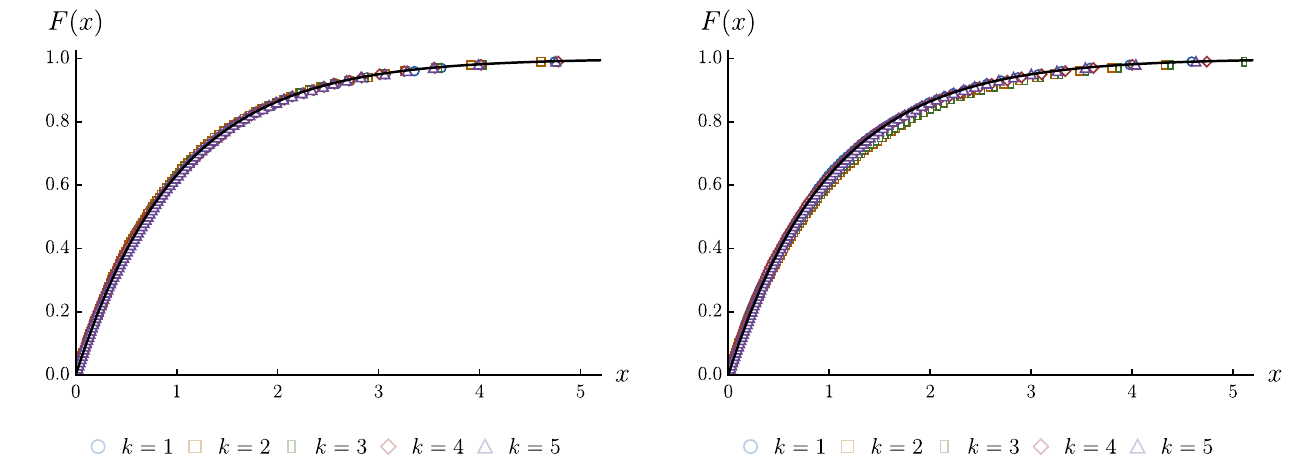}
\caption{Numerical illustrations of the majorisation scaling limits.
Left: empirical distribution function, computed from $10^4$ realisations, of $\mu_{n,k}m\left(\phi_k(\mathsf{\Lambda}_{n,m}-\bot_n)\right)^{n-1}$, with $\mu_{n,k}$ given in~\eqref{eq:maj-limit-one-point-exp}.
Right: empirical distribution function, computed from $10^4$ realisations, of $\nu_{n,k}m(1-\phi_k(\mathsf{V}_{n,m}))^{n-k}$, with $\nu_{n,k}$ given in~\eqref{eq:join-limit-marginals}.
Both simulations are performed for $n=6$ and $m=1000$.
In both panels, the limiting exponential distribution function with rate $1$ is shown as a solid black line.}
\label{fig:meet_join_maj_num}
\end{figure}

We first consider the scaling limit of the majorisation meet as $m\to\infty$.
For this purpose, we denote the cone of nonnegative concave sequences by
\begin{equation}
      \mathcal{C}_{n-1}:=\left\{x=(x_1,\dots,x_{n-1})\in \R^{n-1}_{\geq0}\colon 2 x_k\geq x_{k-1}+x_{k+1}\text{ for all $k\in\llbracket n-1\rrbracket$}\right\}, 
\end{equation}
with the convention $x_0=x_n=0$.
In fact, with this convention the requirement that all coordinates be nonnegative is redundant; see Lemma~\ref{lem:concave-ratios}.
For $s\in \R_{\geq 0}^{n-1}$, define also the shifted orthant 
\begin{equation}
     \mathcal{O}(s)
 :=\left\{x=(x_1,\dots,x_{n-1})\in \R^{n-1}_{\geq0}\colon x_k\geq s_k\text{ for all $k\in\llbracket n-1\rrbracket$}\right\}.
\end{equation}
Throughout, $\operatorname{Vol}$ denotes the Lebesgue measure in the coordinates $(x_1,\ldots,x_{n-1})$.
The following theorem is illustrated in Fig.~\ref{fig:meet_join_maj_num}.

\begin{thm}[Fixed-dimensional scaling limit for the majorisation meet]
\label{thm:maj-conv-distr}
Fix $n\geq 2$.
As $m\to\infty$,
\begin{equation}
m^{\frac{1}{n-1}}\left(\phi_k\left(\mathsf{\Lambda}_{n,m}-\bot_n\right)\right)_{k=1}^{n-1}
\xlongrightarrow{\operatorname{d}}  (\zeta_1,\dots,\zeta_{n-1}),
\end{equation}
where the random vector $(\zeta_1,\dots,\zeta_{n-1})$ is supported on the cone $\mathcal{C}_{n-1}$
and has distribution given by
 \begin{equation}
 \label{eq:maj-volume-formula}
     \P\left(\zeta_1\geq s_1,\ldots, \zeta_{n-1}\geq s_{n-1}\right)=\exp\left(-n!(n-1)!\operatorname{Vol}\left(\mathcal{P}(s)\right)\right)
 \end{equation}
 for all $s=(s_1,\dots,s_{n-1})\in \R_{\geq 0}^{n-1}$, with $\operatorname{Vol}\left(\mathcal{P}(s)\right)$ denoting the volume of the $(n-1)$-dimensional polyhedral region
\begin{equation}
\label{eq:polyhedral_reg}
     \mathcal{P}(s):=\mathcal{C}_{n-1}\setminus \mathcal{O}(s).
\end{equation}
 Moreover, this limiting vector is palindromic in distribution:
\begin{equation}
    (\zeta_1,\dots,\zeta_{n-1}) \stackrel{\mathrm{d}}{=} (\zeta_{n-1},\dots,\zeta_1).
\end{equation}
\end{thm}

\begin{rem}
The limiting law in Theorem~\ref{thm:maj-conv-distr} also has a direct
formulation in the standard coordinates of the simplex.
Since
\begin{equation}
    (\mathsf{\Lambda}_{n,m})_k-\frac1n
=
\phi_k(\mathsf{\Lambda}_{n,m}-\bot_n)
-
\phi_{k-1}(\mathsf{\Lambda}_{n,m}-\bot_n),
\qquad k\in\llbracket n\rrbracket,
\end{equation}
with the convention $\phi_0:=0$, the continuous mapping theorem gives
\begin{equation}
    m^{1/(n-1)}
\left(\mathsf{\Lambda}_{n,m}-\bot_n\right)
\xrightarrow{\operatorname{d}}
\left(
\zeta_1,
\zeta_2-\zeta_1,
\ldots,
\zeta_{n-1}-\zeta_{n-2},
-\zeta_{n-1}
\right).
\end{equation}
Thus, in contrast with the UT-majorisation case (see Remark~\ref{rem:ScalingCoordinatesUT}), no lower-order terms are discarded when passing from partial-sum coordinates to standard coordinates.
\end{rem}

\begin{ex}
Let us compute the volume appearing in the limiting law~\eqref{eq:maj-volume-formula} for $n=2,3,4$.
We may assume that $s\in\mathcal C_{n-1}$.
Indeed, for an arbitrary $s\in\R_{\geq0}^{n-1}$, replacing the augmented sequence $(0,s_1,\ldots,s_{n-1},0)$ by its concave envelope does not change the set of elements of $\mathcal C_{n-1}$ lying componentwise above $s$, and hence does not change the joint survival function in~\eqref{eq:maj-volume-formula}.
A direct computation gives
\begin{align*}
\operatorname{Vol}\left(\mathcal{P}\left(s_1\right)\right) &= s_1, \\
\operatorname{Vol}\left(\mathcal{P}\left(s_1,s_2\right) \right)&=s_1^2-s_1s_2+s_2^2, \\
\operatorname{Vol}\bigl(\mathcal P(s_1,s_2,s_3)\bigr)
&=\frac16\bigl[
-6\underline{s}s_2^2+3\underline{s}^2s_2+4\underline{s}^3
-6(s_2^2-\underline{s}s_2+\underline{s}^2)\overline{s}
+3s_2\overline{s}^2+4s_2^3+2\overline{s}^3
\bigr],
\end{align*}
where, in the last formula,
\[
\underline{s}:=\min(s_1,s_3),
\qquad
\overline{s}:=\max(s_1,s_3).
\]
Already for $n=3$, the mixed term $-s_1s_2$ shows that the joint survival function does not factorise, and hence the coordinates of the limiting vector are dependent.
For $n=4$, the appearance of $\min(s_1,s_3)$ and $\max(s_1,s_3)$ shows that the volume is genuinely
piecewise polynomial rather than polynomial.
The invariance of these formulae under the reverse permutation of the $s_i$ illustrates the palindromic symmetry stated in Theorem~\ref{thm:maj-conv-distr}.
\end{ex}

The following proposition clarifies some properties of $\operatorname{Vol}\bigl(\mathcal P(s)\bigr)$ that emerged in the preceding example.

\begin{prop}[Properties of $\operatorname{Vol}\bigl(\mathcal P(s)\bigr)$]
\label{prop:maj-volume-piecewise}
Let $G=(G_{ij})_{i,j\in\llbracket n-1\rrbracket}$ be the symmetric matrix defined by
\begin{equation}
\label{eq:matrix-G}
G_{ij}
:=
\frac{i\bigl(n-j\bigr)}{n}, \qquad 1\leq i\leq j\leq n-1.
\end{equation}
Then, for every $s\in\mathbb R_{\geq0}^{n-1}$,
\begin{equation}
\label{eq:maj-volume-inclusion-exclusion}
\operatorname{Vol}\bigl(\mathcal P(s)\bigr)
=
\frac{1}{n}
\sum_{\varnothing\neq A\subseteq\llbracket n-1\rrbracket}
(-1)^{|A|+1}
\operatorname{Vol}\bigl(\mathcal{P}_A(s)\bigr),
\end{equation}
where, for each nonempty $A\subseteq\llbracket n-1\rrbracket$, $\mathcal{P}_A(s)$ is the polytope
\[
\mathcal{P}_A(s)
:=
\left\{
z\in\mathbb R_{\geq0}^{n-1}\colon
\sum_{j=1}^{n-1}G_{ij}z_j\leq s_i
\text{ for every }i\in A
\right\}.
\]
In particular, the map $s\mapsto\operatorname{Vol}\bigl(\mathcal P(s)\bigr)$ is continuous and piecewise polynomial on $\mathbb R_{\geq0}^{n-1}$.
\end{prop}

The inclusion--exclusion formula~\eqref{eq:maj-volume-inclusion-exclusion} allows to compute some low-order marginals of the limiting vector, as we next show.

\begin{cor}[Low-order marginals of the limiting vector]
\label{cor:maj-limit-marginals}
Let $(\zeta_1,\ldots,\zeta_{n-1})$ be the limiting vector in Theorem~\ref{thm:maj-conv-distr}.

\begin{enumerate}
\item
For every $k\in\llbracket n-1\rrbracket$,
\begin{equation}
\label{eq:maj-limit-one-point-exp}
\zeta_k^{n-1}
\sim
\operatorname{Exp}(\mu_{n,k}),
\qquad
\mu_{n,k}
:=
\binom{n-1}{k}
\frac{n^{n-1}}
{k^{n-k-1}(n-k)^k}.
\end{equation}

\item
Suppose that $n\geq3$ and, for any $s,t\geq0$, define
\begin{equation}
\label{eq:first-adjacent-concavified-thresholds}
(\tilde{s},\tilde{t})
:=
\begin{cases}
\left(s,\dfrac{n-2}{n-1}s\right),
&
t<\dfrac{n-2}{n-1}s,
\\[3mm]
(s,t),
&
\dfrac{n-2}{n-1}s\leq t\leq2s,
\\[3mm]
\left(\dfrac{t}{2},t\right),
&
t>2s.
\end{cases}
\end{equation}
\end{enumerate}
Then,
\begin{equation}
\label{eq:maj-limit-two-point}
\P(\zeta_1\geq s,\zeta_2\geq t)
=
\P(\zeta_1\geq\tilde{s},\zeta_2\geq\tilde{t})
=
\exp\bigl(-\Psi_n(\tilde{s},\tilde{t})\bigr),
\end{equation}
where
\[
\Psi_n(\tilde{s},\tilde{t})
:=
n^{n-1}\tilde{s}^{n-1}
+
\frac{n}{n-2}
\bigl((n-1)\tilde{t}-(n-2)\tilde{s}\bigr)^{n-1}.
\]
\end{cor}

\begin{rem}
The pair $(\tilde{s},\tilde{t})$ in Corollary~\ref{cor:maj-limit-marginals} consists of the first and second coordinates of the concave envelope of the vector $(0,s,t,0,\dots,0)\in\R_{\geq0}^{n-1}$.
\end{rem}

\begin{rem}
The proof of Corollary~\ref{cor:maj-limit-marginals} is based on the inclusion--exclusion formula~\eqref{eq:maj-volume-inclusion-exclusion}.
The same method yields explicit formulae for the joint distribution of any pair of coordinates $(\zeta_k,\zeta_\ell)$.
The formulae for adjacent coordinates are somewhat simpler, but already become considerably more involved away from the boundary.
We therefore record only the particularly simple case $(\zeta_1,\zeta_2)$; by palindromic symmetry, this also determines the joint distribution of $(\zeta_{n-2},\zeta_{n-1})$.
\end{rem}

We now consider the scaling limit of the majorisation join. 

\begin{thm}[Fixed-dimensional scaling limit for the majorisation join]
\label{thm:maj-join-conv-distr}
Fix $n\geq 2$.
As $m\to\infty$,
\begin{equation}
\label{eq:join-joint-limit}
 \left(
 m^{\frac{1}{n-k}}\bigl(1-\phi_k(\mathsf{V}_{n,m})\bigr)
 \right)_{k=1}^{n-1}
 \xlongrightarrow{\operatorname{d}}
 (\eta_1,\ldots,\eta_{n-1}),
\end{equation}
where $\eta_1,\ldots,\eta_{n-1}$ are independent and
\begin{equation}
\label{eq:join-limit-marginals}
 \eta_k^{n-k}\sim
 \Exp(\nu_{n,k}),
 \qquad
 \nu_{n,k}
 :=\binom{n}{k}\binom{n-1}{k-1},
 \qquad k\in\llbracket n-1\rrbracket.
\end{equation}
Moreover, as $m\to\infty$,
\begin{equation}
\label{eq:concavification-inactive}
 \P\left(
   \phi_k(\mathsf{V}_{n,m})=\max_{1\leq i\leq m}\phi_k(Y_{n,i})
   \quad \text{for all }k\in\llbracket n-1\rrbracket
 \right)\longrightarrow 1.
\end{equation}
\end{thm}

The limit~\eqref{eq:concavification-inactive} tells us that the concavification needed in the construction of the join (see Proposition~\ref{prop:meet-majorisation}) is asymptotically irrelevant in the regime of Theorem~\ref{thm:maj-join-conv-distr}.

\subsection{High-dimensional results}

We now consider the regime in which the dimension $n$ and the sample size $m$ may tend to infinity simultaneously.
We measure the collapse to the extrema in terms of the $p$-norms, defined for $x\in\R^n$ by
\[
\|x\|_p:=\left(\sum_{i=1}^n|x_i|^p\right)^{1/p},
\qquad p\in[1,\infty), \qquad \text{and}\qquad
\|x\|_{\infty}:=\max_{1\leq i\leq n}|x_i|.
\]
The relevant quantities are the total mass outside the last coordinate of the meet and, symmetrically, the total mass outside the first coordinate of the join.
We obtain sharp first-order asymptotics and Gumbel fluctuations for these quantities.
\begin{thm}[High-dimensional scaling limit for UT-majorisation]
\label{thm:UT-m-n-scaling}
Let $(m_n)_{n\geq 2}$ be a sequence of positive integers.
Let $D_n^{\wedge}:=\phi_{n-1}\bigl(\mathsf{\Lambda}_{n,m_n}^{\UT}\bigr)$ and $D_n^{\vee}
:=1-\phi_1\bigl(\mathsf{V}_{n,m_n}^{\UT}\bigr)$.
Set
\[
    a_n^{\UT}:=m_n^{-1/(n-1)}.
\]
As $n\to\infty$, we have
\begin{equation}
\label{eq:UT-m-n-scaling}
    \frac{D_n^{\wedge}}{a_n^{\UT}}
    \xrightarrow{\operatorname{a.s.}}1
    \qquad\text{and}\qquad
    \frac{D_n^{\vee}}{a_n^{\UT}}
    \xrightarrow{\operatorname{a.s.}}1.
\end{equation}
If, in addition, $m_n\to\infty$, then
\begin{equation}
\label{eq:UT-m-n-Gumbel}
    n
    \left(1-\frac{D_n^{\wedge}}
         {a_n^{\UT}}\right)
    \xlongrightarrow{\operatorname{d}}G
    \qquad\text{and}\qquad
    n
    \left(1-\frac{D_n^{\vee}}
         {a_n^{\UT}}\right)
    \xlongrightarrow{\operatorname{d}}G,
\end{equation}
where $G\sim\Gumbel$.
\end{thm}

It turns out that $D_n^{\wedge}$ (resp., $D_n^{\vee}$) in Theorem~\ref{thm:UT-m-n-scaling} measures the distance between the meet and the bottom (resp., between the join and the top); see~\eqref{eq:distance-infty-bottom}-\eqref{eq:distance-1-bottom}-\eqref{eq:distance-top} for more details.
Thus, \eqref{eq:UT-m-n-scaling} yields the following phase diagram determining whether or not the convergence to the extrema takes place.
The relevant parameter is the exponential growth rate of the sample size.
\begin{cor}[High-dimensional behaviour of the UT-majorisation meet and join]
\label{cor:UT-m-n-collapse}
\leavevmode
\begin{enumerate}
\item 
If $(m_n)_{n\geq 2}$ satisfies
\[
\frac{\log m_n}{n}\longrightarrow c\in[0,\infty],
\]
then
\begin{align}
\label{eq:UT-meet-m-n-norm-conv}
\|\mathsf{\Lambda}_{n,m_n}^{\UT}-\botUT_n\|_\infty
= \frac{1}{2} \|\mathsf{\Lambda}_{n,m_n}^{\UT}-\botUT_n\|_1
&\xrightarrow{\operatorname{a.s.}}e^{-c},
\\
\label{eq:UT-join-m-n-norm-conv}
\|\mathsf{V}_{n,m_n}^{\UT}-\topUT_n\|_\infty
= \frac{1}{2} \|\mathsf{V}_{n,m_n}^{\UT}-\topUT_n\|_1
&\xrightarrow{\operatorname{a.s.}}e^{-c},
\end{align}
with the convention $e^{-\infty}=0$.
\item
For every fixed $p\in[1,\infty]$,
\begin{align}
\label{eq:UT-meet-m-n-collapse}
\|\mathsf{\Lambda}_{n,m_n}^{\UT}-\botUT_n\|_p
 \xrightarrow{\operatorname{a.s.}}0
&\quad\Longleftrightarrow\quad
\frac{\log m_n}{n}\longrightarrow\infty,
\\
\label{eq:UT-join-m-n-collapse}
\|\mathsf{V}_{n,m_n}^{\UT}-\topUT_n\|_p
 \xrightarrow{\operatorname{a.s.}}0
&\quad\Longleftrightarrow\quad
\frac{\log m_n}{n}\longrightarrow\infty.
\end{align}
\end{enumerate}
\end{cor}

Corollary~\ref{cor:UT-m-n-collapse} reveals three different regimes.
If $m_n$ grows subexponentially in $n$, then the $\ell^{\infty}$-distances converge to 1, and the $\ell^{1}$-distances converge to 2; thus the meet and join remain asymptotically at maximal distance from the corresponding extrema.
If $m_n=e^{cn+o(n)}$ for some $c\in(0,\infty)$, their $\ell^{\infty}$-distances converge to an intermediate value $e^{-c}$.
Convergence to the bottom and top elements occurs precisely in the superexponential regime $\log m_n/n\to\infty$.

Finally we present the high-dimensional behaviour of the majorisation meet and join.
The behaviour of the meet is described by the following result.
\begin{thm}[High-dimensional behaviour of the majorisation meet]
\label{thm:maj-meet-m-n-collapse}
Let $(m_n)_{n\geq 2}$ be a sequence of positive integers.
Then,
\begin{equation}
\label{eq:maj-meet-m-n-collapse}
    n\bigl\|\mathsf{\Lambda}_{n,m_n}-\bot_n\bigr\|_\infty
    \xrightarrow{\operatorname{a.s.}}0
    \quad\Longleftrightarrow\quad
    \frac{\log m_n}{n}\longrightarrow\infty.
\end{equation}
\end{thm}

For the majorisation join, we obtain a more precise scaling and fluctuations result, analogous to Theorem~\ref{thm:UT-m-n-scaling} for the UT-majorisation join.

\begin{thm}[High-dimensional scaling limit for the majorisation join]
\label{thm:maj-join-m-n-scaling}
Let $(m_n)_{n\geq 2}$ be any sequence of positive integers.
Define $D_n^{\vee}:=1-\phi_1\bigl(\mathsf{V}_{n,m_n}\bigr)$ and set
\begin{equation}
\label{eq:maj-join-an}
    a_n:=(nm_n)^{-1/(n-1)}.
\end{equation}
Then,
\begin{equation}
\label{eq:maj-join-m-n-scaling}
    \frac{D_n^{\vee}}{a_n}
    \xrightarrow{\operatorname{a.s.}}1
\end{equation}
and
\begin{equation}
\label{eq:maj-join-m-n-Gumbel}
    n\left(1-\frac{D_n^{\vee}}{a_n}\right) \xlongrightarrow{\operatorname{d}}G,
\end{equation}
where $G\sim\Gumbel$.
\end{thm}
Notice that no growth assumption on $m_n$ is required for the Gumbel limit here: even for $m_n=1$, the maximum over the $n$ coordinates produces the extreme-value fluctuations.

The almost-sure scaling in Theorem~\ref{thm:maj-join-m-n-scaling} yields the following phase diagram for the convergence to the top.
\begin{cor}[High-dimensional behaviour of the majorisation join]
\label{cor:maj-join-m-n-collapse}
\leavevmode
\begin{enumerate}
\item 
If $(m_n)_{n\geq 2}$ satisfies
\[
    \frac{\log m_n}{n}\longrightarrow c\in[0,\infty],
\]
then
\begin{equation}
    \bigl\|\mathsf{V}_{n,m_n}-\top_n\bigr\|_\infty
    = \frac{1}{2} \bigl\|\mathsf{V}_{n,m_n}-\top_n\bigr\|_1
    \xrightarrow{\operatorname{a.s.}}e^{-c},
\end{equation}
where $e^{-\infty}:=0$.
\item 
For every fixed $p\in[1,\infty]$,
\begin{equation}
\label{eq:maj-join-m-n-collapse}
    \bigl\|\mathsf{V}_{n,m_n}-\top_n\bigr\|_p
    \xrightarrow{\operatorname{a.s.}}0
    \quad\Longleftrightarrow\quad
    \frac{\log m_n}{n}\longrightarrow\infty.
\end{equation}
\end{enumerate}
\end{cor}

\begin{rem}
From~\eqref{eq:maj-meet-m-n-collapse} and~\eqref{eq:maj-join-m-n-collapse}, we see that the same superexponential threshold
governs the convergence of both the majorisation meet and the majorisation join, as for UT-majorisation.
However, in~\eqref{eq:maj-meet-m-n-collapse} the distance $\norm{\mathsf{\Lambda}_{n,m_n}-\bot_n}_{\infty}$ needs to be rescaled by a factor $n$.
The reason is simple to explain: suppose $m_n=1$ for all $n$.
Then, $\mathsf{\Lambda}_{n,1}$ is a point uniformly distributed on $\Delta^\downarrow_{n-1}$, and one can show that, without rescaling,
\begin{equation}
\norm{\mathsf{\Lambda}_{n,1}-\bot_n}_\infty \xrightarrow{\operatorname{a.s.}}  0
\end{equation}
as $n\to\infty$. 
This happens because in the large dimensional limit $n\to\infty$, the uniform measure concentrates around $\bot_n=\left(\frac{1}{n},\ldots,\frac{1}{n}\right)$ in sup-norm.

We also remark that, by~\eqref{eq:maj-meet-m-n-collapse}, the inequality $\|x\|_{\infty}\leq \|x\|_p \leq n^{1/p}\|x\|_{\infty}$ implies 
\begin{align*}
n\norm{\mathsf{\Lambda}_{n,m_n}-\bot_n}_p
\xrightarrow{\operatorname{a.s.}}0
\quad&\Longrightarrow\quad
\frac{\log m_n}{n}\longrightarrow\infty \\
\intertext{and}
\frac{\log m_n}{n}\longrightarrow\infty
\quad&\Longrightarrow\quad
\norm{\mathsf{\Lambda}_{n,m_n}-\bot_n}_p
\xrightarrow{\operatorname{a.s.}}0,
\end{align*}
for any $p\in [1,\infty]$.
However, a faster growing condition on $(m_n)_{n\geq 2}$ might be needed for $n\norm{\mathsf{\Lambda}_{n,m_n}-\bot_n}_1 \xrightarrow{\operatorname{a.s.}}0$.
\end{rem}

\section{Distributional representations and probabilistic structure}
\label{sec:distr_identities}

We first derive general distributional identities for random meets and joins under arbitrary sampling distributions.
We then specialise to the uniform distribution and describe it in partial-sum coordinates, first through explicit integral representations and subsequently through its correspondence with uniform spacings.
Finally, we record the Markov properties of the partial-sum sequences arising in the two settings.

\subsection{General identities for random meets and joins}

The following propositions hold for i.i.d.\ samples from an arbitrary probability distribution on the relevant simplex.

\begin{prop}
Let $X$ be a random point in $\Delta_{n-1}$ with an arbitrary distribution, and let $X_{n,1},X_{n,2},\ldots$ be i.i.d.\ copies of $X$.
Let $\mathsf{\Lambda}_{n,m}^{\UT}:=\bigwedge_{i=1}^m X_{n,i}$ and  $\mathsf{V}_{n,m}^{\UT}:=\bigvee_{i=1}^m X_{n,i}$ in the UT-majorisation lattice. 
Then, for every $(t_1,\dots,t_{n-1})\in\mathbb{R}^{n-1}$:
\begin{align}
\label{eq:distr_meet_UT}
\P\left(\phi_k(\mathsf{\Lambda}_{n,m}^{\UT})\geq t_k\colon k\in \llbracket n-1\rrbracket\right)
&= \P\left(\phi_k(X)\geq t_k \colon k\in \llbracket n-1\rrbracket \right)^m, \\
\label{eq:distr_join_UT}
\P\left(\phi_k(\mathsf{V}_{n,m}^{\UT})\leq t_k\colon k\in \llbracket n-1\rrbracket \right)
&= \P\left(\phi_k(X)\leq t_k\colon k\in \llbracket n-1\rrbracket \right)^m.
\end{align}
\end{prop} 
\begin{proof}
Using Proposition~\ref{prop:meet-UT} and the fact that $X_{n,1},\ldots,X_{n,m}$ are i.i.d., we compute 
\begin{equation*}
\begin{aligned}
&\P\left(\phi_1(\mathsf{\Lambda}_{n,m}^{\UT})\geq t_1,\ldots,\phi_{n-1}(\mathsf{\Lambda}_{n,m}^{\UT})\geq t_{n-1} \right)\\
= \,\,&\P\left(\min_{i\in \llbracket m\rrbracket}\phi_1(X_{n,i})\geq t_1,\min_{i\in \llbracket m\rrbracket}\phi_2(X_{n,i})\geq t_2,\ldots,\min_{i\in \llbracket m\rrbracket}\phi_{n-1}(X_{n,i})\geq t_{n-1}\right)\\
= \,\,&\prod_{i=1}^m \P\left(\phi_1(X_{n,i})\geq t_1,\phi_2(X_{n,i})\geq t_2,\ldots,\phi_{n-1}(X_{n,i})\geq t_{n-1}\right) \\
= \,\,&\P\left(\phi_1(X)\geq t_1,\phi_2(X)\geq t_2,\ldots,\phi_{n-1}(X)\geq t_{n-1}\right)^m,
\end{aligned}
\end{equation*}
which shows~\eqref{eq:distr_meet_UT}.
The proof of~\eqref{eq:distr_join_UT} is similar.
\end{proof}

\begin{prop}
\label{prop:distr-maj}
Let $Y$ be a random point in $\Delta_{n-1}^{\downarrow}$ with an arbitrary distribution, and let $Y_{n,1},Y_{n,2},\ldots$ be i.i.d.\ copies of $Y$.
Let $\mathsf{\Lambda}_{n,m}:=\bigwedge_{i=1}^m Y_{n,i}$ and  $\mathsf{V}_{n,m}:=\bigvee_{i=1}^m Y_{n,i}$ in the majorisation lattice. 
Then, for every $(t_1,\dots,t_{n-1})\in\mathbb{R}^{n-1}$:
\begin{align}
\label{eq:distr_meet_maj}
\P\left(\phi_k(\mathsf{\Lambda}_{n,m})\geq t_k\colon k\in \llbracket n-1\rrbracket\right)
&= \P\left(\phi_k(Y)\geq t_k \colon k\in \llbracket n-1\rrbracket \right)^m, \\
\label{eq:distr_join_maj}
\P\left(\phi_k(\mathsf{V}_{n,m})\leq t_k\colon k\in \llbracket n-1\rrbracket \right)
&\leq \P\left(\phi_k(Y)\leq t_k\colon k\in \llbracket n-1\rrbracket \right)^m.
\end{align}
If $(0,t_1,\ldots,t_{n-1},1)$ is concave, then the two probabilities in~\eqref{eq:distr_join_maj} are equal.
\end{prop} 
\begin{proof}
The identity for the meet follows as in the proof of the preceding proposition.
For the join, set
\[
M:=
\max_{i\in\llbracket m\rrbracket}
\widehat{\phi}(Y_{n,i}),
\qquad
\widehat{t}:=(0,t_1,\ldots,t_{n-1},1),
\]
where the maximum is componentwise.
By Proposition~\ref{prop:meet-majorisation},
\[
\widehat{\phi}(\mathsf{V}_{n,m})=M^*.
\]
Since $M^*\geq M$ and $Y_{n,1},\ldots,Y_{n,m}$ are i.i.d., we have
\[
\begin{aligned}
\P\left(
\phi_k(\mathsf{V}_{n,m})\leq t_k
\colon k\in\llbracket n-1\rrbracket
\right)
&=
\P\left(M^*\leq\widehat{t}\;\right)
\\
&\leq
\P\left(M\leq\widehat{t}\;\right)
\\
&=
\P\left(
\phi_k(Y)\leq t_k
\colon k\in\llbracket n-1\rrbracket
\right)^m,
\end{aligned}
\]
which proves~\eqref{eq:distr_join_maj}.

If $\widehat{t}$ is concave, then $M\leq\widehat{t}$ implies $M^*\leq\widehat{t}$ by the minimality of the concave envelope.
Hence these two events coincide, and \eqref{eq:distr_join_maj} is an equality.
\end{proof}

\subsection{Uniform distributions in partial-sum coordinates}

We derive here integral representations for the uniform distributions on $\Delta_{n-1}$ and $\Delta_{n-1}^{\downarrow}$ in partial-sum
coordinates.

Notice that, in the majorisation lattice, the uniform distribution on the ordered simplex is the natural one induced by classical majorisation on the full simplex.
Indeed, if $X$ is uniformly distributed on $\Delta_{n-1}$, then its decreasing rearrangement $X^\downarrow$ is uniformly distributed on $\Delta_{n-1}^{\downarrow}$.

We start from the following observation.
\begin{lem} \label{lem:bijection}
The images of $\Delta_{n-1}$ and $\Delta_{n-1}^{\downarrow}$ by $\phi$ are
\begin{equation*}
    \phi\left(\Delta_{n-1}\right)=\left\{
u\in\R^n\colon
0\leq u_1\leq\cdots\leq u_{n-1}\leq u_n=1
\right\}
\end{equation*}
and
\begin{equation*}
      \phi(\Delta_{n-1}^{\downarrow})=\left\{u\in\phi\left(\Delta_{n-1}\right)\colon \text{$(0,u_1,\ldots,u_n)$ is concave}  \right\}.
\end{equation*}
\end{lem}
\begin{proof}
By~\eqref{eq:phi-inverse}, the coordinates of $x=\phi^{-1}(u)$ are
\[
x_k=u_k-u_{k-1},
\qquad k\in\llbracket n\rrbracket,
\]
where $u_0:=0$.
Thus $x_k\geq0$ for every $k$ and $\sum_{k=1}^n x_k=1$ if and only if
\[
0=u_0\leq u_1\leq\cdots\leq u_n=1.
\]
Moreover,
\[
x_k-x_{k+1}
=
2u_k-u_{k-1}-u_{k+1}, \qquad k\in \llbracket n-1 \rrbracket.
\]
Consequently, $x_1\geq\cdots\geq x_n$ if and only if
$(u_0,u_1,\ldots,u_n)$ is concave.
\end{proof}

Lemma~\ref{lem:bijection} and the measure-preserving property of $\phi$ yield the following integral representations.
\begin{lem}
\label{lem:unif_UT} Let $X$ be uniformly distributed on $\Delta_{n-1}$.
Then, $(\phi_1(X),\ldots, \phi_{n-1}(X))$ is uniformly distributed over the polytope
\begin{equation}
\label{eq:E_n-1}
    \left\{u\in\R^{n-1}\colon  0\leq u_1\leq\cdots\leq u_{n-1}\leq 1\right\}.
\end{equation}
Equivalently, if we define
\begin{equation}
\label{eq:joint_tail_UT_def}
P^{\UT}_{n-1}(t_1,\ldots,t_{n-1}):=\P\left(\phi_1(X)\geq t_1,\ldots, \phi_{n-1}(X)\geq t_{n-1}\right)
\end{equation}
for $0\leq t_1\leq t_2\leq \cdots \leq t_{n-1}\leq 1$, then
\begin{equation}
    \label{eq:joint_tail_UT}
\begin{aligned}
P^{\UT}_{n-1}(t_1,\ldots,t_{n-1}) 
=(n-1)! \int\limits_{t_{n-1}}^1\diff u_{n-1} \int\limits_{t_{n-2}}^{u_{n-1}}\diff u_{n-2}\cdots \int\limits_{t_{1}}^{u_{2}}\diff u_{1}.
\end{aligned}
\end{equation}
\end{lem}
\begin{proof}
The map
\[
(x_1,\ldots,x_{n-1})
\longmapsto
(\phi_1(x),\ldots,\phi_{n-1}(x))
\]
is linear with Jacobian determinant $1$.
Upon identifying $\Delta_{n-1}$ with its representation in the coordinates $(x_1,\ldots,x_{n-1})$, it has volume $1/(n-1)!$.
Thus, by Lemma~\ref{lem:bijection}, the density of $(\phi_1(X),\ldots, \phi_{n-1}(X))$ is  $(n-1)!$ on the polytope in~\eqref{eq:E_n-1}.
Integrating this density over the region defined by $u_k\geq t_k$ gives~\eqref{eq:joint_tail_UT}.
\end{proof}
The superscript in $P^{\UT}_{n-1}$ indicates that this function arises in UT-majorisation.
The preceding Lemma is illustrated in Fig.~\ref{fig:inf_simplex_UTmaj}.

\begin{ex}
For $n=2,3,4$, the integral in~\eqref{eq:joint_tail_UT} gives
\begin{align*}
P^{\UT}_1(t_1)
&=1-t_1,\\
P^{\UT}_2(t_1,t_2)
&=1-2t_1-t_2^2+2t_1t_2,\\
P^{\UT}_3(t_1,t_2,t_3)
&=1-3t_1-3t_2^2-t_3^3
  +6t_1t_2+3t_1t_3^2+3t_2^2t_3-6t_1t_2t_3.
\end{align*}
By~\eqref{eq:distr_meet_UT}, these expressions determine the joint survival function of the UT-majorisation meet.
More precisely, for $s_1,\dots,s_{n-1}>0$ and large enough $m$,
\[
\P\left(
m\phi_k(\mathsf{\Lambda}_{n,m}^{\UT})^k\geq s_k
\colon k\in\llbracket n-1\rrbracket
\right)
=
P^{\UT}_{n-1}\left(
\frac{s_1}{m},
\left(\frac{s_2}{m}\right)^{1/2},
\ldots,
\left(\frac{s_{n-1}}{m}\right)^{1/(n-1)}
\right)^m.
\]
Consequently, taking exponential limits in the above explicit formulae gives
\begin{align*}
P^{\UT}_1\left(\frac{s_1}{m}\right)^m
&\xrightarrow{m\to\infty}e^{-s_1},\\
P^{\UT}_2\left(
\frac{s_1}{m},
\left(\frac{s_2}{m}\right)^{1/2}
\right)^m
&\xrightarrow{m\to\infty}e^{-2s_1-s_2},\\
P^{\UT}_3\left(
\frac{s_1}{m},
\left(\frac{s_2}{m}\right)^{1/2},
\left(\frac{s_3}{m}\right)^{1/3}
\right)^m
&\xrightarrow{m\to\infty}e^{-3s_1-3s_2-s_3}.
\end{align*}
Thus, the appropriate powers of the limiting coordinates are independent and exponentially distributed, in agreement with Theorem~\ref{thm:UT-conv-distr}.
\end{ex}

\begin{lem}
\label{lem:concave-ratios}
Let $(u_0,\ldots,u_n)$ be a concave real sequence with $u_0,u_n\geq0$.
Then
\be
\label{eq:concave-ratios}
\frac{u_1}{1}\geq\frac{u_2}{2}\geq\cdots\geq\frac{u_{n-1}}{n-1}\geq\frac{u_n}{n},
\qquad
\frac{u_0}{n}\leq\frac{u_1}{n-1}\leq\cdots\leq\frac{u_{n-1}}{1}.
\ee
\end{lem}
\begin{proof}
Set $d_j:=u_j-u_{j-1}$ for $j=1,\ldots,n$.
By concavity, $d_1\geq\cdots\geq d_n$.
Using $u_0\geq 0$, for $k\in\llbracket n-1\rrbracket$ we thus have
\[
\frac{u_k}{k}
\geq \frac{u_k-u_0}{k}
=
\frac{d_1+\cdots+d_k}{k}
\geq d_{k+1}
=
u_{k+1}-u_k,
\]
and hence $u_{k+1}/(k+1)\leq u_k/k$, proving the first chain.
For the second chain, it suffices to apply the first chain to the reversed sequence $(u_n,\dots,u_0)$, which is still concave and, by hypothesis, satisfies $u_n\geq 0$.
\end{proof}

\begin{lem}
\label{lem:unif_maj}
Let $Y$ be uniformly distributed on $\Delta_{n-1}^{\downarrow}$.
Then $(\phi_1(Y),\ldots,\phi_{n-1}(Y))$ is uniformly distributed over the polytope
\[
\mathcal{A}_{n-1}
:=\left\{
u\in\R^{n-1}\colon
0\leq u_1\leq\cdots\leq u_{n-1}\leq1,
\;\;
2u_k\geq u_{k-1}+u_{k+1},
\;\; k\in\llbracket n-1\rrbracket
\right\},
\]
with the usual convention $u_0:=0$, $u_n:=1$. Equivalently, this polytope can be written as
\begin{equation}
\label{eq:polytope-maj-recursive}
\mathcal{A}_{n-1}=
\left\{
u\in\R^{n-1}\colon
\frac{n-1}{n}\leq u_{n-1}\leq1,
\,
\frac{k}{k+1}u_{k+1}
\leq u_k
\leq2u_{k+1}-u_{k+2},
\, k\in\llbracket n-2\rrbracket
\right\}.
\end{equation}

For $0\leq t_1\leq\cdots\leq t_{n-1}\leq1$, define
\begin{align*}
P_{n-1}(t_1,\ldots,t_{n-1})
&:=
\P\left(
\phi_k(Y)\geq t_k
\colon k\in\llbracket n-1\rrbracket
\right),\\
Q_{n-1}(t_1,\ldots,t_{n-1})
&:=
\P\left(
\phi_k(Y)\leq t_k
\colon k\in\llbracket n-1\rrbracket
\right).
\end{align*}
Then
\begin{multline}
\label{eq:joint_tail_maj}
P_{n-1}(t_1,\ldots,t_{n-1})
=n!(n-1)!\times\\
\int_{\max(\frac{n-1}{n},t_{n-1})}^{1}\!\!\diff u_{n-1}
\int_{\max(\frac{n-2}{n-1}u_{n-1},t_{n-2})}^{2u_{n-1}-1}
   \!\!\diff u_{n-2}
\cdots
\int_{\max(\frac{2}{3}u_3,t_2)}^{2u_3-u_4}\!\!\diff u_2
\int_{\max(\frac{1}{2}u_2,t_1)}^{2u_2-u_3}\!\!\diff u_1,
\end{multline}
whereas
\begin{multline}
\label{eq:joint_cdf_maj}
Q_{n-1}(t_1,\ldots,t_{n-1})
=n!(n-1)!\times\\
\int_{\frac{n-1}{n}}^{t_{n-1}}\!\!\diff u_{n-1}
\int_{\frac{n-2}{n-1}u_{n-1}}^{\min(2u_{n-1}-1,t_{n-2})}
   \!\!\diff u_{n-2}
\cdots
\int_{\frac{2}{3}u_3}^{\min(2u_3-u_4,t_2)}\!\!\diff u_2
\int_{\frac{1}{2}u_2}^{\min(2u_2-u_3,t_1)}\!\!\diff u_1.
\end{multline}
In both formulae, an integral $\int_a^b$ with $a>b$ is understood to be zero.
\end{lem}

\begin{proof}
The same change of variables as before has Jacobian determinant $1$.
Under the same identification with the coordinates $(x_1,\ldots,x_{n-1})$, coordinate permutations partition $\Delta_{n-1}$, up to a set of measure zero, into $n!$ chambers of equal volume.
Since $\Delta_{n-1}^{\downarrow}$ is one of these chambers, it has volume $\frac{1}{n!(n-1)!}$.
Thus, by Lemma~\ref{lem:bijection}, the image of its uniform distribution under $\phi$ has density $n!(n-1)!$ on the polytope $\mathcal{A}_{n-1}$.

To verify~\eqref{eq:polytope-maj-recursive}, suppose first that $u$ belongs to $\mathcal{A}_{n-1}$.
Then, $(0,u_1,\dots,u_{n-1},1)$ is concave, so the lower bounds in~\eqref{eq:polytope-maj-recursive} follow from Lemma~\ref{lem:concave-ratios}; the upper bounds follow immediately from the definition of $\mathcal{A}_{n-1}$.

Conversely, suppose that $u$ belongs to the polytope on the right-hand side of~\eqref{eq:polytope-maj-recursive}.
We have $u_{n-1}\leq u_n:=1$ and, if we suppose by backward induction that $u_{k+1}\leq u_{k+2}$ for some $k\in \llbracket n-2\rrbracket$, then
\[
u_k \leq 2u_{k+1} -u_{k+2} \leq u_{k+1}.
\]
The lower bounds imply $u_k\geq0$ for every $k$.
Moreover, the upper bounds are precisely the concavity inequalities at indices $2,\ldots,n-1$, whereas the lower bound $u_1\geq u_2/2$ is the concavity inequality at index $1$.
Thus $(0,u_1,\ldots,u_{n-1},1)$ is nondecreasing and concave, proving that $u\in \mathcal{A}_{n-1}$.

Finally, intersecting the recursive bounds in~\eqref{eq:polytope-maj-recursive} with $u_k\geq t_k$ (respectively, with $u_k\leq t_k$) gives the iterated integral in~\eqref{eq:joint_tail_maj} (respectively, in~\eqref{eq:joint_cdf_maj}).
\end{proof}

For an illustration of Lemma~\ref{lem:unif_maj}, see Fig.~\ref{fig:inf_simplex_maj}.

The functions $P_{n-1}$ and $Q_{n-1}$ determine, respectively, the joint upper and lower tails of the partial sums of a uniformly distributed point in $\Delta_{n-1}^{\downarrow}$.
By Proposition~\ref{prop:distr-maj}, they are related to the majorisation meet and join by
\begin{align*}
\P\left(
\phi_k(\mathsf{\Lambda}_{n,m})\geq t_k
\colon k\in\llbracket n-1\rrbracket
\right)
&=
P_{n-1}(t_1,\ldots,t_{n-1})^m, \\
\P\left(
\phi_k(\mathsf{V}_{n,m})\leq t_k
\colon k\in\llbracket n-1\rrbracket
\right)
&\leq
Q_{n-1}(t_1,\ldots,t_{n-1})^m.
\end{align*}
If the augmented sequence $(0,t_1,\ldots,t_{n-1},1)$ is concave, then the latter inequality is an equality.
Thus $P_{n-1}$ always gives the exact distribution of the meet, whereas $Q_{n-1}$ gives the exact distribution of the join at concave thresholds and an upper bound at general thresholds.
The distinction is caused precisely by the concavification involved in constructing the majorisation join.

\subsection{Spacings of uniform points}
\label{sec:spacings}

The explicit integral formulae~\eqref{eq:joint_tail_UT}, \eqref{eq:joint_tail_maj}, \eqref{eq:joint_cdf_maj} for the joint distributions of uniform points on the simplex and the ordered simplex in partial-sum coordinates are not especially well suited to asymptotic analysis.
We therefore use the classical correspondence between uniform points on the simplex and uniform spacings obtained by placing $n-1$ independent uniform points in the unit interval; see, for example,~\cite{Renyi53,Pyke1965}.
For further results on ordered uniform spacings, see~\cite{Bairamov2010}.

Suppose in the following that $(U_i)_{i\geq1}$ is a sequence of i.i.d.\ random variables uniformly distributed on the interval $[0,1]$.
Let $n\geq 2$ and consider the first $n-1$ variables $U_1,\ldots, U_{n-1}$.
Let
\[
0=:U_{0:n-1}\leq U_{1:n-1}\leq U_{2:n-1}\leq \cdots\leq U_{n-1:n-1}\leq U_{n:n-1}:=1
\]
and
\[
S_{i,n}:=U_{i:n-1}-U_{i-1:n-1}, \qquad i=1,\ldots, n,
\]
be, respectively, the order statistics and the $n$ spacings determined by these $n-1$ points.
Arranging the spacings $S_{1,n},\ldots, S_{n,n}$ in decreasing order, we  obtain ordered uniform spacings $S_{1,n}^{\downarrow}\geq S_{2,n}^{\downarrow}\geq\cdots\geq S_{n,n}^{\downarrow}$. We denote by $S_n=(S_{1,n},\ldots, S_{n,n})$ and $S^{\downarrow}_n=(S_{1,n}^{\downarrow},\ldots, S_{n,n}^{\downarrow})$ the vectors of spacings and ordered spacings, respectively.

\begin{prop}
\label{prop:UT-uniform}
Let $X=(X_1,\ldots,X_n)$ be uniformly distributed on $\Delta_{n-1}$.
Then
\begin{align}
     \label{eq:UT-uniform}
     (X_1,\ldots,X_n)&\stackrel{\mathrm{d}}{=}(S_{1,n},\ldots, S_{n,n}), \\
     \label{eq:UT-uniform-partial}
     (\phi_1(X),\ldots,\phi_{n-1}(X))&\stackrel{\mathrm{d}}{=}\left(U_{1:n-1},\ldots, U_{n-1:n-1}\right).
\end{align}
\end{prop}
\begin{proof}
By Lemma~\ref{lem:unif_UT}, the vector $(\phi_1(X),\ldots,\phi_{n-1}(X))$ has density $(n-1)!$ on the region~\eqref{eq:E_n-1}.
The order statistics $(U_{1:n-1},\ldots,U_{n-1:n-1})$ have the same density on that region, and~\eqref{eq:UT-uniform-partial} follows.
Since $\phi_k(S_n)=U_{k:n-1}$ for each $k\in\llbracket n-1\rrbracket$ and the first $n-1$ partial sums uniquely determine a point of $\Delta_{n-1}$, we conclude that~\eqref{eq:UT-uniform} also holds.
\end{proof}

Thus, for $0\leq t_1\leq \cdots\leq t_{n-1}\leq 1$, we have
\[
     P^{\UT}_{n-1}(t_1,\dots,t_{n-1})
     =\P\left(U_{1:n-1}\geq t_1,\ldots, U_{n-1:n-1}\geq t_{n-1}\right),
\]
where $P^{\UT}_{n-1}(t_1,\dots,t_{n-1})$ is defined in~\eqref{eq:joint_tail_UT_def}.

\begin{lem}
\label{lem:joint-upper-tails-spacings}
For every $I\subseteq\llbracket n\rrbracket$ and $(t_i)_{i\in I}\in\mathbb R_{\geq0}^{I}$,
\[
\P\left(S_{i,n}\geq t_i\colon i\in I\right)
=
\left(1-\sum_{i\in I}t_i\right)^{n-1}
\mathbbm{1}_{\left\{\sum_{i\in I}t_i\leq1\right\}}.
\]
\end{lem}
\begin{proof}
Set $T:=\sum_{i\in I}t_i$.
If $T>1$, the event is empty.
Suppose therefore that $T\leq1$.
The corresponding subset of the simplex is
\[
\left\{
s\in\Delta_{n-1}\colon s_i\geq t_i\text{ for every }i\in I
\right\}
=
\left(\sum_{i\in I}t_i e_i\right)
+
(1-T)\Delta_{n-1},
\]
where $(e_1,\dots,e_n)$ is the canonical basis of $\R^n$.
Indeed, subtracting $t_i$ from the $i$-th coordinate for each $i\in I$ leaves nonnegative coordinates with total mass $1-T$.
Translation preserves $(n-1)$-dimensional volume, while dilation by $1-T$ multiplies it by $(1-T)^{n-1}$.
Since $S_n$ is uniformly distributed on $\Delta_{n-1}$ by Proposition~\ref{prop:UT-uniform}, the required probability is therefore
\[
\frac{\operatorname{Vol}\bigl((1-T)\Delta_{n-1}\bigr)}
     {\operatorname{Vol}(\Delta_{n-1})}
=
(1-T)^{n-1},
\]
as desired.
\end{proof}

There is a useful recursion for the joint survival function of uniform order statistics, discovered by Bolshev and subsequently published in~\cite{Kotelnikova83}; see also~\cite[\S~9.3]{Shorak09} and~\cite[Proposition~5.7]{Cunden25}.
\begin{prop}[Bolshev's recursion]
For any $0\leq t_1\leq \cdots\leq t_{n-1}\leq 1$,
    \begin{equation}
    \label{eq:Bolshev}
    \begin{aligned}
P^{\UT}_{n-1}(t_1,\dots,t_{n-1}) 
= 1- \sum_{k=1}^{n-1} \binom{n-1}{k} t_k^k P^{\UT}_{n-1-k}(t_{k+1},\dots, t_{n-1}),
\end{aligned}
\end{equation}
with $P^{\UT}_0:=1$. 
\end{prop}
The one-point marginal distributions can be easily computed; see e.g.~\cite{OrderStatistics}. 
\begin{lem}
\label{lem:marginal_phik_UT}
Let $X$ be uniformly distributed on $\Delta_{n-1}$.
Then, for all $t\in [0,1]$,
\[
     \P\left(\phi_k(X)\leq t\right)=\P\left(U_{k:n-1}\leq t\right)=
     \sum_{i=k}^{n-1}\binom{n-1}{i}t^i(1-t)^{n-1-i}, \qquad k\in\llbracket n-1\rrbracket.
\]
Equivalently, $\phi_k(X)$ has the $\operatorname{Beta}(k,n-k)$ distribution.
In particular, for $k=1$ and $k=n-1$,
\begin{equation}
    \label{eq:marginal_1_n-1_UT}
    \P\left(\phi_1(X)\leq 1-t\right)=   \P\left(\phi_{n-1}(X)\geq t\right)=
    1-t^{n-1}, \qquad t\in [0,1].
\end{equation}
\end{lem}

Let us now consider analogous properties on the ordered simplex.

\begin{prop} 
\label{prop:maj-uniform}
Let $Y=(Y_1,\ldots,Y_n)$ be uniformly distributed on $\Delta_{n-1}^{\downarrow}$.
Then
\begin{equation}
     (Y_1,\ldots,Y_n)\stackrel{\mathrm{d}}{=}(S_{1,n}^{\downarrow},\ldots, S_{n,n}^{\downarrow}),
\end{equation}
and therefore
\begin{equation}
     (\phi_1(Y),\ldots,\phi_{n-1}(Y))\stackrel{\mathrm{d}}{=}(\phi_1(S^{\downarrow}_n),\ldots,\phi_{n-1}(S^{\downarrow}_n)).
\end{equation}
\end{prop}
\begin{proof}
By Proposition~\ref{prop:UT-uniform}, $S_n$ is uniformly distributed on $\Delta_{n-1}$.
Its decreasing rearrangement is therefore uniformly distributed on the permutation chamber $\Delta_{n-1}^{\downarrow}$.
\end{proof}

We are not aware of an explicit recursion for the joint distribution of the ordered spacings in partial-sum coordinates, in the same spirit as~\eqref{eq:Bolshev}.
The one-point marginals, however, are known explicitly: they coincide with the distributions of the sums of the $k$ largest uniform spacings, computed classically by Mauldon~\cite{Mauldon1951}; see also~\cite{ShtembariCaldwell2020} for a recent derivation.
\begin{lem} 
\label{lem:marginal_phik_maj}
Let $Y$ be uniformly distributed on $\Delta_{n-1}^{\downarrow}$.
Then, for all $t\in [0,1]$ and $k\in\llbracket n-1\rrbracket$,
\begin{equation}
\label{eq:marginal_phik_maj}
\begin{split}
     &\P\left(\phi_k(Y)\leq t\right)=\P\left(S_{1,n}^{\downarrow}+\cdots+S_{k,n}^{\downarrow}\leq t\right)=\\
     &\frac{n! }{k!k^{n-k-1}}\sum_{i=k}^{n-1}\frac{(-1)^{i-k} (n-i)^{n-k-1}}{ (n+k-i) (n-i-1)!(i-k)!}\left(\frac{t(n+k-i)-k}{n-i}\right)^{n-1} \mathbbm{1}_{\left[\frac{k}{n+k-i},1\right]}(t). 
\end{split}
\end{equation}
\end{lem}

\begin{cor}
\label{cor:endpoint-marginals-maj}
Let $Y$ be uniformly distributed on $\Delta_{n-1}^{\downarrow}$.
Then
\begin{align}
\label{eq:marginal_1_n-1_maj}
\P\left(\phi_1(Y)\leq1-t\right)
&=
1-nt^{n-1},
&&0<t\leq\frac12,
\\
\label{eq:marginal_n-1_lower_maj}
\P\left(
\phi_{n-1}(Y)\geq\frac{n-1}{n}+t
\right)
&=
1-(nt)^{n-1},
&&0<t\leq\frac1n,
\\
\label{eq:marginal_1_lower_maj}
\P\left(
\phi_1(Y)\geq\frac1n+t
\right)
&=
1-(nt)^{n-1},
&&0<t\leq\frac{1}{n(n-1)}.
\end{align}
\end{cor}
\begin{proof}
The last two identities follow directly from Lemma~\ref{lem:marginal_phik_maj}.
For $k=n-1$, the formula in that lemma reduces to
\be
\label{eq:marginal_n-1_upper_maj}
\P\left(\phi_{n-1}(Y)\leq u\right)
=
\left(nu-(n-1)\right)^{n-1},
\qquad
\frac{n-1}{n}\leq u\leq1.
\ee
Setting $u=(n-1)/n+t$ and taking complements gives~\eqref{eq:marginal_n-1_lower_maj}.

For $k=1$ and $1/n\leq u\leq1/(n-1)$, only the first term in the sum in Lemma~\ref{lem:marginal_phik_maj} is nonzero, and hence
\[
\P\left(\phi_1(Y)\leq u\right)
=
(nu-1)^{n-1}.
\]
Setting $u=1/n+t$ and taking complements gives~\eqref{eq:marginal_1_lower_maj}.

The first identity~\eqref{eq:marginal_1_n-1_maj} can also be deduced from the lemma, but a direct argument is easier.
Let $X=(X_1,\ldots,X_n)$ be uniformly distributed on $\Delta_{n-1}$.
Since $Y\stackrel{\mathrm d}=X^\downarrow$,
\[
\phi_1(Y)\stackrel{\mathrm d}=\max_{j\in\llbracket n\rrbracket}X_j.
\]
For $0<t\leq1/2$, the events $\{X_j>1-t\}$, $j\in\llbracket n\rrbracket$, are pairwise disjoint and, by~\eqref{eq:marginal_1_n-1_UT}, have probability $t^{n-1}$.
Therefore
\[
\P\left(\phi_1(Y)\leq1-t\right)
=
1-\sum_{j=1}^n\P(X_j>1-t)
=
1-nt^{n-1},
\]
which proves~\eqref{eq:marginal_1_n-1_maj}.
\end{proof}

\subsection{Markov properties}
\label{sec:Markov}

We conclude by recording the Markov structure of the partial-sum sequences under the two uniform distributions.

\begin{prop} 
\label{prop:UT-Markov}
Let $X$ be uniformly distributed on $\Delta_{n-1}$.
Then the sequence
\[
(\phi_n(X),\phi_{n-1}(X), \dots, \phi_1(X))
\]
is a time-inhomogeneous Markov chain with initial condition $\phi_n(X)=1$ and transition densities
\begin{equation}
\label{eq:transition_UT}
    f_{\phi_{k}|\phi_{k+1}}(x|y)=\frac{k}{y}\left(\frac{x}{y}\right)^{k-1}\mathbbm{1}_{[0,y]}(x), 
    \qquad k\in\llbracket n-1\rrbracket.
\end{equation}
\end{prop}
\begin{proof}
The identity $\phi_n(X)=1$ follows from the definition of $\phi$.
By Proposition~\ref{prop:UT-uniform},
\[
(\phi_1(X),\ldots,\phi_{n-1}(X))
\stackrel{\mathrm{d}}{=}
(U_{1:n-1},\ldots,U_{n-1:n-1}).
\]
Conditionally on $U_{k+1:n-1}=y$, the first $k$ order statistics have the same distribution as the order statistics of $k$ independent uniform random variables on $[0,y]$, independently of the subsequent order statistics.
This already shows that the conditional distribution
of $U_{k:n-1}$ given \( U_{k+1:n-1},\ldots,U_{n-1:n-1} \) depends only on $U_{k+1:n-1}$, and hence proves the Markov property.
Moreover, $U_{k:n-1}$ is the largest of these $k$ uniform random variables: it lies below some threshold $x$ if and only if all of them are below the same threshold.
Therefore, for $0\leq x\leq y$,
\[
\mathbb P\!\left(
U_{k:n-1}\leq x
\,\middle|\,
U_{k+1:n-1}=y
\right)
=
\left(\frac{x}{y}\right)^k.
\]
Differentiating with respect to $x$ gives the conditional density~\eqref{eq:transition_UT}.
\end{proof}

\begin{prop}
\label{prop:maj-2-Markov}
Let $Y$ be uniformly distributed on $\Delta_{n-1}^{\downarrow}$.
Then the sequence
\[
\bigl(\phi_n(Y),\phi_{n-1}(Y),\ldots,\phi_1(Y)\bigr)
\]
is a time-inhomogeneous second-order Markov chain.
Its two initial states are such that $\phi_n(Y)=1$ almost surely and $\phi_{n-1}(Y)$ has density
\begin{equation}
\label{eq:initial_maj}
f_{\phi_{n-1}}(y)=(n-1)n^{n-1}\left(y-\frac{n-1}{n}\right)^{n-2}\mathbbm{1}_{\left[\frac{n-1}{n},1\right]}(y).
\end{equation}
For $k\in\llbracket n-2\rrbracket$, its transition densities are
\begin{equation}
\label{eq:transition_maj}
f_{\phi_k\mid\phi_{k+1},\phi_{k+2}}(x\mid y,z)
=
\frac{k(k+1)\bigl((k+1)x-ky\bigr)^{k-1}}
     {\bigl((k+2)y-(k+1)z\bigr)^k}
\mathbbm{1}_{\left[\frac{k}{k+1}y,\,2y-z\right]}(x),
\qquad k\in\llbracket n-2\rrbracket.
\end{equation}
\end{prop}

\begin{proof}
The identity $\phi_n(Y)=1$ follows from the definition of $\phi$.
The density~\eqref{eq:initial_maj} of $\phi_{n-1}(Y)$ follows from differentiating~\eqref{eq:marginal_n-1_upper_maj}.

Recall that, by Lemma~\ref{lem:unif_maj}, the vector $(\phi_1(Y),\ldots,\phi_{n-1}(Y))$ has constant density on $\mathcal A_{n-1}$.
Fix $k\in\llbracket n-2\rrbracket$ and condition on $(\phi_{k+1}(Y),\ldots,\phi_n(Y))$, writing $\phi_{k+1}(Y)=y$ and $\phi_{k+2}(Y)=z$.
If $\phi_k(Y)=x$, then $Y_{k+1}=y-x$ and $Y_{k+2}=z-y$, so the ordering $Y_{k+1}\geq Y_{k+2}$ requires $x\leq2y-z$.
The first $k$ coordinates are decreasing, have sum $x$, and are bounded below by $y-x$.
After subtracting $y-x$ from each of them, they form a decreasing nonnegative $k$-tuple with sum $(k+1)x-ky$.
Such tuples exist if and only if $x\geq\frac{k}{k+1}y$, and their volume in the first $k-1$ coordinates is proportional to $\bigl((k+1)x-ky\bigr)^{k-1}$.

Since the joint density on $\mathcal A_{n-1}$ is constant, the conditional density of $\phi_k(Y)$ given $(\phi_{k+1}(Y),\ldots,\phi_n(Y))$ is therefore proportional to $\bigl((k+1)x-ky\bigr)^{k-1}$ on $\bigl[\frac{k}{k+1}y,2y-z\bigr]$ and depends only on $y$ and $z$, proving the second-order Markov property.
Finally,
\[
\int_{\frac{k}{k+1}y}^{2y-z}\bigl((k+1)x-ky\bigr)^{k-1}\,\mathrm{d}x
=
\frac{\bigl((k+2)y-(k+1)z\bigr)^k}{k(k+1)},
\]
which gives~\eqref{eq:transition_maj}.
\end{proof}

The different Markov orders reflect the geometry of the corresponding partial-sum polytopes.
In the UT-majorisation case, the admissible range of $u_k$, given the subsequent coordinates, is determined solely by $u_{k+1}$ through
\[
0\leq u_k\leq u_{k+1}.
\]
In the majorisation case, the recursive bounds take the form
\[
\frac{k}{k+1}u_{k+1}\leq u_k\leq 2u_{k+1}-u_{k+2},
\]
so the conditional distribution of $u_k$ generally depends on both $u_{k+1}$ and $u_{k+2}$.
Thus the passage from a first-order to a second-order Markov structure is a direct probabilistic manifestation of the concavity constraint imposed by majorisation.

\section{Proofs}
\label{sec:proofs}

We will repeatedly use the inequality
 \begin{align}
\label{eq:bounds_norms} 
    \|x\|_p\leq \|x\|_r,\quad \text{for } x\in\R^{n}, \quad 1\leq r\leq p\leq \infty.
\end{align}

\subsection{Proofs of the main results on UT-majorisation}
\label{sec:proofs_UTmaj}
Recall that, in the UT-majorisation relation, $\botUT_{n}=(0,\ldots,0,1)$ and $\topUT_{n}=(1,\ldots,0,0)$.
Note that, for all $x\in\Delta_{n-1}$, since $x_j\geq 0$ for all $j$ and $1-x_n=x_1+\cdots+x_{n-1}$, we have
\begin{align}
\label{eq:distance-infty-bottom}
    \norm{x-\botUT_n}_\infty
    &=\max\{x_1,\dots,x_{n-1},1-x_n\}=1-x_n=\phi_{n-1}(x), \\
\label{eq:distance-1-bottom}
    \| x-\botUT_n\|_{1}
    &=x_1+\cdots+x_{n-1}+(1-x_n)=2\phi_{n-1}(x).
\end{align}
Analogously
\begin{equation}
\label{eq:distance-top}
    \norm{x-\topUT_n}_\infty =1-\phi_{1}(x), \qquad \norm{x-\topUT_n}_1=2(1-\phi_1(x)).
\end{equation}
Therefore the inequalities~\eqref{eq:bounds_norms} imply 
\be
\label{eq:p-ineq-bottom}
\phi_{n-1}(x)
\leq \norm{x-\botUT_n}_p
\leq 
 2\phi_{n-1}(x)
\leq 2(1-\phi_{1}(x)).
\ee
and
\be
\label{eq:p-ineq-top}
1-\phi_{1}(x)
\leq \norm{x-\topUT_n}_p
\leq 2(1-\phi_{1}(x)).
\ee

\subsubsection{Fixed-dimensional results}
\begin{proof}[Proof of Theorem~\ref{thm:UT-conv-as}]
From~\eqref{eq:distance-1-bottom} and~\eqref{eq:distance-top}, we have
\begin{align*}
\norm{\mathsf{\Lambda}_{n,m}^{\UT}-\botUT_n}_1 \xrightarrow{{\text{a.s.}}} 0 \quad&\iff\quad \phi_{n-1}(\mathsf{\Lambda}_{n,m}^{\UT})\xrightarrow{{\text{a.s.}}}  0 \\
\intertext{and}
\norm{\mathsf{V}_{n,m}^{\UT}-\topUT_n}_1\xrightarrow{{\text{a.s.}}} 0 \quad&\iff\quad \phi_{1}(\mathsf{V}_{n,m}^{\UT}) \xrightarrow{{\text{a.s.}}}  1.
\end{align*}
From~\eqref{eq:distr_meet_UT}, \eqref{eq:distr_join_UT}, and~\eqref{eq:marginal_1_n-1_UT}, we obtain, for every $t\in(0,1)$,
\[
\P\left(\phi_{n-1}(\mathsf{\Lambda}_{n,m}^{\UT})\geq t\right)
=
\P\left(\phi_1(\mathsf{V}_{n,m}^{\UT})\leq 1-t\right)
=
\left(1-t^{n-1}\right)^m.
\]
Therefore, as $m\to\infty$, $\phi_{n-1}(\mathsf{\Lambda}_{n,m}^{\UT})\to 0$ and $\phi_{1}(\mathsf{V}_{n,m}^{\UT})\to 1$ in probability.
Since $\phi_{n-1}(\mathsf{\Lambda}_{n,m}^{\UT})$ (resp., $\phi_{1}(\mathsf{V}_{n,m}^{\UT})$)  is a nonincreasing (resp., nondecreasing) sequence in $m$, convergence in
probability implies almost sure convergence.
\end{proof}

\begin{proof}[Proof of Theorem~\ref{thm:UT-conv-distr}]
Recall the definition of $P^{\UT}_{n-1}$ in~\eqref{eq:joint_tail_UT_def}.
Assuming without loss of generality that $s_1,\dots,s_{n-1}>0$ and $m$ is large enough, we want to show that 
\begin{equation}
\label{eq:UTmeet_asymptotics}
\lim_{m\to\infty}P^{\UT}_{n-1}\left(\frac{s_1}{m}, \left(\frac{s_2}{m}\right)^{\frac{1}{2}}, \dots, \left(\frac{s_{n-1}}{m}\right)^{\frac{1}{n-1}}\right)^{m}=\prod_{k=1}^{n-1}e^{-\binom{n-1}{k}s_k}.
\end{equation}
Using Bolshev's recursion~\eqref{eq:Bolshev} and the limit $P^{\UT}_{j}(t_1,\dots,t_j)\to 1$ as $t_1,\dots,t_j\to 0$, we compute
\[
\begin{aligned}
&\quad \log\left( P^{\UT}_{n-1}\left(\frac{s_1}{m}, \left(\frac{s_2}{m}\right)^{\frac{1}{2}}, \dots, \left(\frac{s_{n-1}}{m}\right)^{\frac{1}{n-1}}\right)^m\right) \\
&= m\log \left(1- \sum_{k=1}^{n-1} \binom{n-1}{k} \frac{s_k}{m} P^{\UT}_{n-1-k}\left(\left(\frac{s_{k+1}}{m}\right)^{\frac{1}{k+1}}, \dots, \left(\frac{s_{n-1}}{m}\right)^{\frac{1}{n-1}}\right) \right) \\
&\xrightarrow{m\to\infty} -\sum_{k=1}^{n-1} \binom{n-1}{k} s_k.
\end{aligned}
\]
This proves~\eqref{eq:UTmeet_asymptotics}.

For the join, consider the reflection \eqref{eq:Reflection}. 
This transformation preserves the uniform distribution on $\Delta_{n-1}$ and satisfies
\[
\phi_k(Rx)=1-\phi_{n-k}(x),
\qquad k\in\llbracket n-1\rrbracket.
\]
Consequently, if $\widetilde X_{n,i}:=R(X_{n,i})$, then the $\widetilde X_{n,i}$ are again i.i.d.\ uniform points and
\[
\begin{aligned}
\phi_k\left(\bigwedge_{i=1}^m\widetilde X_{n,i}\right)
=\min_{i\in\llbracket m\rrbracket}\phi_k(\widetilde X_{n,i})
=1-\max_{i\in\llbracket m\rrbracket}\phi_{n-k}(X_{n,i})
=1-\phi_{n-k}(\mathsf{V}_{n,m}^{\UT}).
\end{aligned}
\]
Therefore,
\[
\left(\phi_k(\mathsf{\Lambda}_{n,m}^{\UT})\right)_{k=1}^{n-1}
\stackrel{\mathrm d}{=}
\left(1-\phi_{n-k}(\mathsf{V}_{n,m}^{\UT})\right)_{k=1}^{n-1},
\]
and the result for the join follows from the result for the meet.
\end{proof}

\subsubsection{High-dimensional results}

\begin{proof}[Proof of Theorem~\ref{thm:UT-m-n-scaling}]
By~\eqref{eq:distr_meet_UT}, \eqref{eq:distr_join_UT}, and~\eqref{eq:marginal_1_n-1_UT}, for every $t\in[0,1]$,
\[
\P\left(D_n^{\wedge}\geq t\right)
=
\P\left(D_n^{\vee}\geq t\right)
=
\left(1-t^{n-1}\right)^{m_n}.
\]
Denote by $D_n$ any of the two random variables $D_n^{\wedge}$ and $D_n^{\vee}$.
Then, the distribution of the random variable
\[
B_n:=\left(D_n\right)^{n-1}
\]
is given by
\begin{equation}
    \label{eq:UT-beta}
\P\left(B_n\geq t\right)
=
(1-t)^{m_n},
\qquad t\in[0,1].
\end{equation}

To prove~\eqref{eq:UT-m-n-scaling}, notice first that the almost sure convergence
\[
\frac{D_n}{a_n^{\UT}}
=
\left(m_nB_n\right)^{1/(n-1)}
\xrightarrow{\operatorname{a.s.}}1
\]
is equivalent to
\begin{equation}
    \label{eq:UT-as-conv-proof}
\frac{\log(m_nB_n)}{n-1}\xrightarrow{\operatorname{a.s.}}0.
\end{equation}
Fix $\epsilon>0$.
Using~\eqref{eq:UT-beta} and the elementary inequality $1-(1-u)^m\leq mu$, valid for $u\in[0,1]$, we obtain
\[
\P\left(\frac{\log(m_nB_n)}{n-1}<-\epsilon\right)
= \P\left(B_n<\frac{e^{-\epsilon(n-1)}}{m_n}\right)
=1-\left(1-\frac{e^{-\epsilon(n-1)}}{m_n}\right)^{m_n}
\leq e^{-\epsilon(n-1)}.
\]
On the other hand, using~\eqref{eq:UT-beta} and the inequality $1-x\leq e^{-x}$ for $x\geq0$ we obtain
\[
\P\left(\frac{\log(m_nB_n)}{n-1}>\epsilon\right)
= \P\left(B_n>\frac{e^{\epsilon(n-1)}}{m_n}\right)
=
\left(1-\frac{e^{\epsilon(n-1)}}{m_n}\right)^{m_n}
\leq \exp\left(-e^{\epsilon(n-1)}\right)
\]
if $e^{\epsilon(n-1)}<m_n$, whereas this probability is zero if $e^{\epsilon(n-1)}\geq m_n$.
Therefore,
\[
\sum_{n=2}^{\infty}
\P\left(
\left|\log(m_nB_n)\right|>\epsilon(n-1)
\right)
<\infty.
\]
Applying the Borel--Cantelli lemma yields~\eqref{eq:UT-as-conv-proof}.

Suppose now that $m_n\to\infty$.
For every $x\geq0$ and all sufficiently large $n$,
\[
\P\left(m_nB_n>x\right)
=
\left(1-\frac{x}{m_n}\right)^{m_n}
\longrightarrow e^{-x}.
\]
Therefore,
\[
m_nB_n\xrightarrow{\operatorname{d}}E,
\qquad
E\sim\operatorname{Exp}(1).
\]
The continuous mapping theorem thus gives
\[
G_n:=-\log\left(m_nB_n\right)
\xrightarrow{\operatorname{d}}-\log E=:G.
\]
Moreover,
\[
\P\left(G\leq x\right)
=
\P\left(E\geq e^{-x}\right)
=
\exp\left(-e^{-x}\right),
\qquad x\in\R,
\]
and hence $G\sim\Gumbel$.

Finally,
\[
\frac{D_n}{a_n^{\UT}}
=
\exp\left(-\frac{G_n}{n-1}\right).
\]
Since $G_n$ is bounded in probability, a Taylor expansion gives
\[
n\left(1-\frac{D_n}{a_n^{\UT}}\right)
=
n\left(1-\exp\left(-\frac{G_n}{n-1}\right)\right)
=
G_n+o(1),
\]
where $o(1)$ denotes a random variable converging to zero in probability.
Slutsky's theorem therefore yields
\[
n\left(1-\frac{D_n}{a_n^{\UT}}\right)
\xrightarrow{\operatorname{d}}G,
\]
which proves~\eqref{eq:UT-m-n-Gumbel}.
\end{proof}

\begin{proof}[Proof of Corollary~\ref{cor:UT-m-n-collapse}]
Since
\[
\log a_n^{\UT}
=
-\frac{\log m_n}{n-1},
\]
we have
\begin{equation}
\label{eq:UT-m-n-critical-scaling}
\frac{\log m_n}{n} \to c
\qquad\Longleftrightarrow\qquad
a_n^{\UT} \to e^{-c}
\end{equation}
for any $c\in[0,\infty]$, with the convention $e^{-\infty}:=0$.

Using the notation $D_n^{\wedge}:= \phi_{n-1}(\mathsf{\Lambda}_{n,m_n}^{\UT})$ and $D_n^\vee:=1-\phi_1(\mathsf{V}_{n,m_n}^{\UT})$ as in Theorem~\ref{thm:UT-m-n-scaling}, it follows from~\eqref{eq:bounds_norms}--\eqref{eq:p-ineq-top} that, for every $p\in[1,\infty]$,
\begin{align}
\label{eq:UT-meet-p-inequalities}
D_n^\wedge
&=
\bigl\|\mathsf{\Lambda}_{n,m_n}^{\UT}-\botUT_n\bigr\|_\infty
\leq
\bigl\|\mathsf{\Lambda}_{n,m_n}^{\UT}-\botUT_n\bigr\|_p
\leq
\bigl\|\mathsf{\Lambda}_{n,m_n}^{\UT}-\botUT_n\bigr\|_1
=
2D_n^\wedge \\
\label{eq:UT-join-p-inequalities}
D_n^\vee
&=
\bigl\|\mathsf{V}_{n,m_n}^{\UT}-\topUT_n\bigr\|_\infty
\leq
\bigl\|\mathsf{V}_{n,m_n}^{\UT}-\topUT_n\bigr\|_p
\leq
\bigl\|\mathsf{V}_{n,m_n}^{\UT}-\topUT_n\bigr\|_1
=
2D_n^\vee.
\end{align}

If $\log m_n/n\to c\in[0,\infty]$, then $a_n^{\UT}\to e^{-c}$, and hence the almost-sure convergence~\eqref{eq:UT-m-n-scaling} gives $D_n^\wedge\xrightarrow{\operatorname{a.s.}} e^{-c}$ and $D_n^\vee\xrightarrow{\operatorname{a.s.}} e^{-c}$, proving~\eqref{eq:UT-meet-m-n-norm-conv}--\eqref{eq:UT-join-m-n-norm-conv}.

Moreover, using~\eqref{eq:UT-meet-p-inequalities},
\eqref{eq:UT-m-n-scaling}, and~\eqref{eq:UT-m-n-critical-scaling}, we have
\[
\bigl\|\mathsf{\Lambda}_{n,m_n}^{\UT}-\botUT_n\bigr\|_p
\xrightarrow{\operatorname{a.s.}}0
\quad\Longleftrightarrow\quad
D_n^\wedge\xrightarrow{\operatorname{a.s.}}0
\quad\Longleftrightarrow\quad
a_n^{\UT}\to0
\quad\Longleftrightarrow\quad
\frac{\log m_n}{n}\to\infty,
\]
for every $p\in[1,\infty]$.
This proves~\eqref{eq:UT-meet-m-n-collapse}.
The proof of~\eqref{eq:UT-join-m-n-collapse} is similar.
\end{proof}

\subsection{Proofs of the main results on majorisation}
\label{sec:proofs_maj}

Recall that, in the majorisation relation, $\bot_n=\left(\frac1n,\ldots,\frac1n\right)$ and $\top_n=(1,0,\ldots,0)$.
For every $x\in\Delta_{n-1}^{\downarrow}$ and every $p\in[1,\infty]$, the monotonicity $x_1\geq x_2\geq\cdots\geq x_n$ implies
\begin{equation}
\label{eq:meet_maj_ineq}
\begin{split}
\max\left\{
x_1-\frac1n,\frac1n-x_n
\right\}
=
\norm{x-\bot_n}_\infty
&\leq
\norm{x-\bot_n}_p \\
&\leq
\norm{x-\bot_n}_1
=
2\norm{\phi(x)-\phi(\bot_n)}_\infty.
\end{split}
\end{equation}
To verify the last identity, let $r:=\max\left\{
j\in\llbracket n\rrbracket\colon x_j\geq\frac{1}{n}
\right\}$ and note that
\[
\norm{x-\bot_n}_1
=\sum_{j=1}^r\left(x_j-\frac1n\right) + \sum_{j=r+1}^n\left(\frac1n-x_j\right)
=2\left(\phi_r(x)-\frac rn\right)
=2\norm{\phi(x)-\phi(\bot_n)}_\infty.
\]
Indeed, the successive differences of the sequence $\bigl(\phi_k(x)-k/n\bigr)_{k=0}^n$ are $x_k-1/n$, so this sequence increases up to $k=r$ and decreases thereafter.

Moreover, as for UT-majorisation, we have
\begin{equation}
\label{eq:join_maj_ineq}
1-\phi_1(x)
=
\norm{x-\top_n}_\infty
\leq
\norm{x-\top_n}_p
\leq
\norm{x-\top_n}_1
=
2\left(1-\phi_1(x)\right).
\end{equation}

\subsubsection{Fixed-dimensional results}

\begin{proof}[Proof of Theorem~\ref{thm:maj-conv-as}]
We first prove the convergence of the meet to $\bot_n$.
Since
\[
(\mathsf{\Lambda}_{n,m})_1
\geq
(\mathsf{\Lambda}_{n,m})_2
\geq
\cdots
\geq
(\mathsf{\Lambda}_{n,m})_n,
\]
it is enough to show that $(\mathsf{\Lambda}_{n,m})_1$ and $(\mathsf{\Lambda}_{n,m})_n$ both converge almost surely to $1/n$.

For the first coordinate, using~\eqref{eq:distr_meet_maj} and~\eqref{eq:marginal_1_lower_maj}, we obtain
\begin{align*}
\P\left(
\abs{(\mathsf{\Lambda}_{n,m})_1-\frac1n}\geq\epsilon
\right)
=\P\left(
\phi_1(\mathsf{\Lambda}_{n,m})\geq\frac1n+\epsilon
\right)
=\left(
1-(n\epsilon)^{n-1}
\right)^m
\xrightarrow{m\to\infty}0
\end{align*}
for $0<\epsilon\leq 1/(n(n-1))$.
Hence $(\mathsf{\Lambda}_{n,m})_1\to1/n$ in probability.
Since the sequence $\bigl((\mathsf{\Lambda}_{n,m})_1\bigr)_{m\geq1}$ is nonincreasing, this convergence holds almost surely.

For the last coordinate, using~\eqref{eq:distr_meet_maj} and~\eqref{eq:marginal_n-1_lower_maj}, we similarly obtain
\begin{align*}
\P\left(
\abs{(\mathsf{\Lambda}_{n,m})_n-\frac1n}\geq\epsilon
\right)
=\P\left(
\phi_{n-1}(\mathsf{\Lambda}_{n,m})
\geq\frac{n-1}{n}+\epsilon
\right)
=\left(
1-(n\epsilon)^{n-1}
\right)^m
\xrightarrow{m\to\infty}0.
\end{align*}
for $0<\epsilon\leq1/n$.
Thus $(\mathsf{\Lambda}_{n,m})_n\to1/n$ in probability.
Since the sequence $\bigl((\mathsf{\Lambda}_{n,m})_n\bigr)_{m\geq1}$ is nondecreasing, this convergence also holds almost surely.

We now prove the almost sure convergence of the join to $\top_n$.
By~\eqref{eq:join_maj_ineq},
\[
\norm{\mathsf{V}_{n,m}-\top_n}_1
\xrightarrow{\mathrm{a.s.}}0
\quad\Longleftrightarrow\quad
\phi_1(\mathsf{V}_{n,m})
\xrightarrow{\mathrm{a.s.}}1,
\]
so it suffices to show the convergence on the right.
Fix $\epsilon\in(0,1/2]$, and consider the probability
\[
\P\left(\abs{\phi_1(\mathsf{V}_{n,m})-1}\geq \epsilon\right)
=\P\left(\phi_1(\mathsf{V}_{n,m})\leq 1-\epsilon\right)
=\P\left(\phi_k(\mathsf{V}_{n,m})\leq t_k\colon k\in \llbracket n-1\rrbracket \right),
\]
where $t_1:=1-\epsilon$ and $t_k=1$ for $2\leq k\leq n-1$.
The augmented sequence $(0,t_1,\dots,t_{n-1},1)$ is readily verified to be concave.
Therefore, the equality case of~\eqref{eq:distr_join_maj}, together with~\eqref{eq:marginal_1_n-1_maj}, gives
\begin{align*}
\P\left(\abs{\phi_1(\mathsf{V}_{n,m})-1}\geq \epsilon\right)
=\P\left(\phi_1(X)\leq 1-\epsilon\right)^m
=\left(
1-n\epsilon^{n-1}
\right)^m
\xrightarrow{m\to\infty}0,
\end{align*}
where $X$ is uniformly distributed on $\Delta^{\downarrow}_{n-1}$.
Hence $\phi_1(\mathsf{V}_{n,m})\to1$ in probability.
Since the sequence $\bigl(\phi_1(\mathsf{V}_{n,m})\bigr)_{m\geq1}$ is nondecreasing, the convergence holds almost surely.
\end{proof}

\begin{proof}[Proof of Theorem~\ref{thm:maj-conv-distr}]
Let $Y$ be uniformly distributed on $\Delta_{n-1}^{\downarrow}$ and $Z=(Z_1,\ldots,Z_{n-1})$ be defined by
\[
Z_k:=\phi_k(Y)-\frac{k}{n},
\qquad
k\in\llbracket n-1\rrbracket.
\]
By Lemma~\ref{lem:unif_maj}, the random vector $Z$ is uniformly distributed on
\[
\widetilde{\mathcal{A}}_{n-1}
:=
\left\{
u\in\mathbb R_{\geq0}^{n-1}\colon
u_{n-1}\leq\frac1n,
\quad
2u_k\geq u_{k-1}+u_{k+1}
\text{ for every }k\in\llbracket n-1\rrbracket
\right\},
\]
where $u_0=u_n=0$.
Moreover,
\[
\operatorname{Vol}\bigl(\widetilde{\mathcal A}_{n-1}\bigr)
=\operatorname{Vol}\bigl(\Delta_{n-1}^{\downarrow}\bigr)
=\frac{1}{n!(n-1)!}.
\]

Since the augmented sequences $(0,u_1,\ldots,u_{n-1},0)$ and $(0,u_{n-1},\ldots,u_1,0)$ are concave, Lemma~\ref{lem:concave-ratios} yields
\begin{equation}
\label{eq:conditionOnus}
\frac{u_1}{1}\geq\frac{u_2}{2}\geq\cdots\geq\frac{u_{n-1}}{n-1},
\qquad
\frac{u_1}{n-1}\leq\frac{u_2}{n-2}\leq\cdots\leq\frac{u_{n-1}}{1}.
\end{equation}

Fix $s=(s_1,\ldots,s_{n-1})\in\mathbb R_{\geq0}^{n-1}$, and set
\[
\mathcal{O}_m(s)
:=
\mathcal O\left(m^{-\frac{1}{n-1}}s\right).
\]
We have $\operatorname{Vol}(\lambda D)=\lambda^{n-1}\operatorname{Vol}(D)$ for all measurable sets $D\subseteq\R^{n-1}$ and $\lambda>0$.
Therefore, \eqref{eq:distr_meet_maj} gives
\begin{equation}
\label{eq:proof-maj-meet-survival}
\begin{aligned}
&\P\left(
m^{\frac{1}{n-1}}
\left(
\phi_k(\mathsf{\Lambda}_{n,m})-\frac{k}{n}
\right)
\geq s_k
\colon k\in\llbracket n-1\rrbracket
\right)
=
\P\left(
Z_k\geq m^{-\frac{1}{n-1}}s_k
\colon k\in\llbracket n-1\rrbracket
\right)^m \\
&=\left(
1-\frac{\operatorname{Vol}\left(\widetilde{\mathcal{A}}_{n-1}\setminus \mathcal{O}_m(s)\right)}
{\operatorname{Vol}\left(\widetilde{\mathcal{A}}_{n-1}\right)}
\right)^m
=\left(
1-\frac{n!(n-1)!}{m}
\operatorname{Vol}\left(
m^{\frac{1}{n-1}}\left(\widetilde{\mathcal{A}}_{n-1}\setminus \mathcal{O}_m(s)\right)
\right)
\right)^m.
\end{aligned}
\end{equation}

We claim that, for all sufficiently large $m$,
\begin{equation}
\label{eq:claim_sets}
m^{\frac{1}{n-1}}\left(\widetilde{\mathcal{A}}_{n-1}\setminus \mathcal{O}_m(s)\right)
=
\mathcal P(s).
\end{equation}
We refer the reader to Fig.~\ref{fig:sets_claim}, where the claim is illustrated in dimension $n=3$.
\begin{figure}
\centering
\includegraphics[width=.7\textwidth]{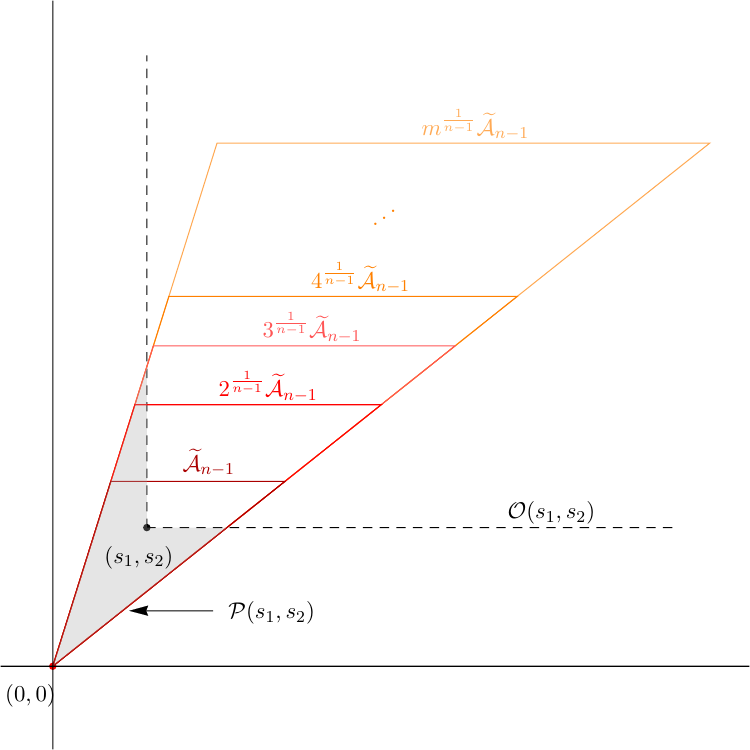}
        \caption{
        The increasing sets $m^{\frac{1}{n-1}}\left(\widetilde{\mathcal{A}}_{n-1}\setminus \mathcal{O}_m(s)\right)$ stabilise to $\mathcal{P}(s)$ (in grey) for $m$ sufficiently large; see~\eqref{eq:claim_sets}.
        An illustration for $n=3$.}
\label{fig:sets_claim}
\end{figure}
Since
\[
m^{\frac{1}{n-1}}\mathcal{O}_m(s)=\mathcal O(s)
\qquad \text{and} \qquad
m^{\frac{1}{n-1}}\widetilde{\mathcal{A}}_{n-1}
=
\left\{
u\in\mathcal C_{n-1}\colon
u_{n-1}\leq\frac{m^{\frac{1}{n-1}}}{n}
\right\},
\]
we have
\[
m^{\frac{1}{n-1}}\left(\widetilde{\mathcal{A}}_{n-1}\setminus \mathcal{O}_m(s)\right)
\subseteq
\mathcal P(s)
\]
for all $m\geq 1$.
For the reverse inclusion, let $u\in\mathcal P(s)$.
Then $u_j<s_j$ for some $j\in\llbracket n-1\rrbracket$.
By the first chain of inequalities in~\eqref{eq:conditionOnus},
\[
u_{n-1}
\leq
\frac{n-1}{j}u_j
<
\frac{n-1}{j}s_j
\leq
C(s):= \max_{i\in\llbracket n-1\rrbracket}
\frac{n-1}{i}s_i.
\]
Choose $m_0\in\mathbb N$ such that $m_0^{\frac{1}{n-1}}/n\geq C(s)$.
Then, for every $m\geq m_0$, $u_{n-1}\leq m^{\frac{1}{n-1}}/n$, so that $u\in m^{\frac{1}{n-1}}\widetilde{\mathcal{A}}_{n-1}$.
Since $u\notin\mathcal O(s)=m^{\frac{1}{n-1}}\mathcal{O}_m(s)$, we conclude that
\[
u\in m^{\frac{1}{n-1}}\left(\widetilde{\mathcal{A}}_{n-1}\setminus \mathcal{O}_m(s)\right),
\]
which proves~\eqref{eq:claim_sets}.

Notice also that $\mathcal P(s)$ is bounded.
Indeed, the preceding argument bounds $u_{n-1}$, while the second chain of inequalities in~\eqref{eq:conditionOnus} gives $u_j\leq(n-j)u_{n-1}$ for all $j\in\llbracket n-1\rrbracket$.
Thus $\operatorname{Vol}\bigl(\mathcal P(s)\bigr)<\infty$.

Substituting~\eqref{eq:claim_sets} into~\eqref{eq:proof-maj-meet-survival}, we obtain, for all $m\geq m_0$,
\begin{align*}
\P\left(
m^{\frac{1}{n-1}}
\left(
\phi_k(\mathsf{\Lambda}_{n,m})-\frac{k}{n}
\right)
\geq s_k
\colon k\in\llbracket n-1\rrbracket
\right)
&=
\left(
1-\frac{n!(n-1)!}{m}
\operatorname{Vol}\bigl(\mathcal P(s)\bigr)
\right)^m \\
&\xrightarrow{m\to\infty}
\exp\left(
-n!(n-1)!\operatorname{Vol}\bigl(\mathcal P(s)\bigr)
\right).
\end{align*}

For $s=ae_k$ (where $e_k$ is the $k$-th vector of the canonical basis of $\R^{n-1}$), the homogeneity of $\mathcal P$ gives
\[
\operatorname{Vol}\bigl(\mathcal P(ae_k)\bigr)
=
a^{n-1}\operatorname{Vol}\bigl(\mathcal P(e_k)\bigr),
\]
where $0<\operatorname{Vol}\bigl(\mathcal P(e_k)\bigr)<\infty$.
Thus the limiting one-point survival functions tend to zero at infinity.
The sequence of random vectors is therefore tight, and the convergence of the multivariate survival functions thus yields the asserted convergence in distribution; in particular, the right-hand side of~\eqref{eq:maj-volume-formula} is the survival function of a probability measure on $\R^{n-1}$.
Since each rescaled vector belongs to the closed cone $\mathcal C_{n-1}$, the limiting vector is supported on $\mathcal C_{n-1}$.

It remains to prove the palindromic symmetry.
The function $R$ defined in \eqref{eq:Reflection} is a measure-preserving involution on $\R^{n-1}$, maps $\mathcal C_{n-1}$ onto itself and satisfies $R\bigl(\mathcal O(s)\bigr)=\mathcal O(Rs)$.
Hence $R\bigl(\mathcal P(s)\bigr)=\mathcal P(Rs)$, which implies $\operatorname{Vol}\bigl(\mathcal P(s)\bigr)=\operatorname{Vol}\bigl(\mathcal P(Rs)\bigr)$.
The survival function~\eqref{eq:maj-volume-formula} is thus invariant under reversal of the coordinates, proving the required palindromic symmetry.
\end{proof}

\begin{proof}[Proof of Proposition~\ref{prop:maj-volume-piecewise}]
Consider the linear map
\[
x\longmapsto z,
\qquad
z_i:=2x_i-x_{i-1}-x_{i+1},
\quad i\in\llbracket n-1\rrbracket,
\]
where $x_0=x_n=0$.
Its matrix is the $(n-1)\times (n-1)$ discrete Laplacian
\[
L
=
\begin{pmatrix}
2 & -1 & 0 & \cdots & 0\\
-1 & 2 & -1 & \ddots & \vdots\\
0 & -1 & 2 & \ddots & 0\\
\vdots & \ddots & \ddots & \ddots & -1\\
0 & \cdots & 0 & -1 & 2
\end{pmatrix}.
\]
The concavity inequalities defining $\mathcal C_{n-1}$ are precisely the inequalities $z_i\geq0$.
Thus, this map sends $\mathcal C_{n-1}$ bijectively onto $\mathbb R_{\geq0}^{n-1}$.

The inverse of $L$ is the Green matrix $G$ defined in~\eqref{eq:matrix-G}.
Expanding along the last row, it is readily shown by induction on $n$ that $\det L=n$.
Consequently, the change of variables $x=Gz$ has Jacobian $1/n$, and
\[
\diff x=\frac1n\diff z.
\]

By definition,
\[
\mathcal P(s)
=
\bigcup_{i=1}^{n-1}
\left\{
x\in\mathcal C_{n-1}\colon x_i<s_i
\right\}.
\]
Under the change of variables $x=Gz$, the inclusion--exclusion principle gives
\[
\operatorname{Vol}\bigl(\mathcal P(s)\bigr)
=
\frac1n
\sum_{\varnothing\neq A\subseteq\llbracket n-1\rrbracket}
(-1)^{|A|+1}
\operatorname{Vol}
\left\{
z\in\mathbb R_{\geq0}^{n-1}\colon
\sum_{j=1}^{n-1}G_{ij}z_j<s_i
\text{ for every }i\in A
\right\}.
\]
The volume of each set in the sum equals that of its closure $\mathcal P_A(s)$, obtained by replacing $<$ with $\leq$.
This proves~\eqref{eq:maj-volume-inclusion-exclusion}.

It remains to prove that the map $s\mapsto\operatorname{Vol}\bigl(\mathcal P(s)\bigr)$ is continuous and piecewise polynomial.
Since~\eqref{eq:maj-volume-inclusion-exclusion} is a finite sum, it suffices to prove the properties for the map $s\mapsto\operatorname{Vol}\bigl(\mathcal P_A(s)\bigr)$, for any fixed nonempty set $A\subseteq\llbracket n-1\rrbracket$.

To prove that $s\mapsto \operatorname{Vol} (\mathcal P_A(s))$ is piecewise polynomial, observe that the coefficients in the inequalities defining $\mathcal P_A(s)$ are independent of $s$, while their right-hand sides depend linearly on $s$.
Every vertex is obtained by imposing equality in $n-1$ linearly independent constraints selected from $z_j\geq0$, $j\in\llbracket n-1\rrbracket$, and $\sum_{j=1}^{n-1}G_{ij}z_j\leq s_i$, $i\in A$.
Each choice of such inequalities produces a linear system with a fixed coefficient matrix and a right-hand side depending linearly on $s$.
Whenever the system is nonsingular, its solution is therefore a linear function of $s$.
Whether a candidate vertex belongs to $\mathcal P_A(s)$ is determined by substituting it into the remaining inequalities.
Since the candidate vertex depends linearly on $s$, these conditions are linear inequalities in $s$.
The set $\mathbb R_{\geq0}^{n-1}$ of possible parameters $s$ can consequently be partitioned into finitely many polyhedral regions on each of which the vertices and the face structure of $\mathcal P_A(s)$ remain unchanged.

On each such region, $\mathcal P_A(s)$ admits a triangulation whose simplices are determined by fixed collections of its vertices.
The vertices depend linearly on $s$, so the determinant formula for the volume of a simplex shows that the volume of each simplex is a polynomial in $s$.
Therefore, $s\mapsto\operatorname{Vol}\bigl(\mathcal P_A(s)\bigr)$ is piecewise polynomial.

Continuity of $s\mapsto \operatorname{Vol}\bigl(\mathcal P_A(s)\bigr)$ follows from the dominated convergence theorem.
Indeed, if $s^{(r)}\to s$ in $\R_{\geq 0}^{n-1}$, then the polytopes $\mathcal P_A(s^{(r)})$ lie in a common bounded set, while their indicator functions converge almost everywhere to that of $\mathcal P_A(s)$, since their boundaries are contained in finitely many hyperplanes.
Thus $\operatorname{Vol}\bigl(\mathcal P_A(s^{(r)})\bigr)\to\operatorname{Vol}\bigl(\mathcal P_A(s)\bigr)$.
\end{proof}

\begin{proof}[Proof of Corollary~\ref{cor:maj-limit-marginals}]
Let $e_k$ denote the $k$-th vector of the canonical basis of $\mathbb R^{n-1}$.
By Proposition~\ref{prop:maj-volume-piecewise},
\begin{equation}
\label{eq:proof-one-point-volume}
\operatorname{Vol}\bigl(\mathcal P(se_k)\bigr)
=
\frac{1}{n}
\operatorname{Vol}
\left\{
z\in\mathbb R_{\geq0}^{n-1}
\colon
\sum_{j=1}^{n-1}G_{kj}z_j\leq s
\right\}.
\end{equation}

For $s>0$, make the change of variables
\[
y_j:=\frac{G_{kj}}{s}z_j,
\qquad
j\in\llbracket n-1\rrbracket.
\]
The resulting region is the standard simplex $\{y\in \R_{\geq 0}^{n-1}\colon \sum_{j=1}^{n-1} y_j<1\}$, whose volume is $1/(n-1)!$, while the Jacobian of the change of variables is $s^{n-1}\prod_{j=1}^{n-1}G_{kj}^{-1}$.
Consequently,
\[
\operatorname{Vol}\bigl(\mathcal P(se_k)\bigr)
=
\frac{s^{n-1}}
{n(n-1)!\prod_{j=1}^{n-1}G_{kj}}.
\]
The same formula is immediate for $s=0$.
Since
\begin{equation*}
\prod_{j=1}^{n-1}G_{kj}
=
\frac{
k!(n-k-1)!\,
k^{n-k-1}(n-k)^k
}
{n^{n-1}},
\end{equation*}
we obtain
\begin{equation}
\label{eq:one-point-marginal-volume}
n!(n-1)!
\operatorname{Vol}\bigl(\mathcal P(se_k)\bigr)
=
\mu_{n,k}s^{n-1}.
\end{equation}
The survival function in Theorem~\ref{thm:maj-conv-distr} therefore gives
\begin{equation}
\P(\zeta_k\geq s)
=
\exp\bigl(-\mu_{n,k}s^{n-1}\bigr),
\end{equation}
which proves~\eqref{eq:maj-limit-one-point-exp}.

We now consider the pair $(\zeta_1,\zeta_2)$.
Suppose first that
\begin{equation}
\label{eq:proof-first-adjacent-compatible}
\frac{n-2}{n-1}s\leq t\leq2s.
\end{equation}
Set
\begin{align*}
V_1(s)
&:=
\operatorname{Vol}
\left\{
z\in\mathbb R_{\geq0}^{n-1}
\colon
(Gz)_1\leq s
\right\},
\\
V_2(t)
&:=
\operatorname{Vol}
\left\{
z\in\mathbb R_{\geq0}^{n-1}
\colon
(Gz)_2\leq t
\right\},
\\
V_{12}(s,t)
&:=
\operatorname{Vol}
\left\{
z\in\mathbb R_{\geq0}^{n-1}
\colon
(Gz)_1\leq s,\ (Gz)_2\leq t
\right\}.
\end{align*}
The inclusion--exclusion formula in Proposition~\ref{prop:maj-volume-piecewise} gives
\begin{equation}
\label{eq:proof-first-pair-inclusion-exclusion}
\operatorname{Vol}\bigl(\mathcal P(se_1+te_2)\bigr)
=
\frac1n
\bigl(
V_1(s)+V_2(t)-V_{12}(s,t)
\bigr).
\end{equation}
Indeed, only the terms indexed by $\{1\}$, $\{2\}$, $\{1,2\}$ survive in inclusion--exclusion: any term involving another index has zero full-dimensional volume because the corresponding threshold is zero and all entries of $G$ are positive.
By~\eqref{eq:proof-one-point-volume} and~\eqref{eq:one-point-marginal-volume}, we know that
\begin{equation}
\label{eq:proof-first-coordinate-volume}
\frac{1}{n}V_1(s)
=\operatorname{Vol}\bigl(\mathcal P(se_1)\bigr)
=\frac{1}{n!(n-1)!} \mu_{n,1} s^{n-1}
=\frac{n^{n-1}s^{n-1}}{n!(n-1)!}.
\end{equation}
We next compute $V_2(t)-V_{12}(s,t)$ as the volume of
\[
\mathcal R(s,t)
:=
\left\{
z\in\mathbb R_{\geq0}^{n-1}
\colon
(Gz)_1\geq s,\ (Gz)_2\leq t
\right\}.
\]
This is simpler than computing $V_{12}(s,t)$ directly.
Set
\[
w_j:=G_{1j}z_j=\frac{n-j}{n}z_j,
\qquad
j\in\llbracket n-1\rrbracket,
\]
so that
\[
(Gz)_1=\sum_{j=1}^{n-1}w_j, \qquad
(Gz)_2=\frac{n-2}{n-1}w_1+2\sum_{j=2}^{n-1}w_j.
\]
The Jacobian of this change of variables is
\[
\frac{\partial z}{\partial w}
=\prod_{j=1}^{n-1}\frac{n}{n-j}
=\frac{n^{n-1}}{(n-1)!}.
\]
Setting
\[
a:=\sum_{j=1}^{n-1}w_j
\qquad\text{and}\qquad
b:=\sum_{j=2}^{n-1}w_j,
\]
the inequalities defining $\mathcal R(s,t)$ are equivalent to
\[
s\leq a\leq \frac{n-1}{n-2}t,
\qquad
0\leq b\leq \frac{n-1}{n}t-\frac{n-2}{n}a.
\]
The assumption $\frac{n-2}{n-1}s\leq t$ from~\eqref{eq:proof-first-adjacent-compatible} ensures that the region $\mathcal R(s,t)$ is nonempty.
On the other hand, the assumption $t\leq2s$ ensures that the upper bound on $b$ is at most $a$, so the condition $a-b=w_1\geq0$ is automatic.
The variables $(a,b,w_2,\ldots,w_{n-2})$ form a system of coordinates, since
\[
w_1=a-b,
\qquad
w_{n-1}=b-\sum_{j=2}^{n-2}w_j,
\]
and the corresponding change of variables has Jacobian determinant $1$.
For fixed $a$ and $b$, the admissible values of the free coordinates are therefore
\[
\left\{
(w_2,\ldots,w_{n-2})\in\mathbb R_{\geq0}^{n-3}
\colon
\sum_{j=2}^{n-2}w_j\leq b
\right\},
\]
whose volume is $b^{n-3}/(n-3)!$.
Hence
\begin{align*}
\operatorname{Vol}\bigl(\mathcal R(s,t)\bigr)
&=
\frac{n^{n-1}}
{(n-1)!}
\int_s^{\frac{n-1}{n-2}t} \diff a
\int_0^{\frac{n-1}{n}t-\frac{n-2}{n}a} \diff b 
\frac{b^{n-3}}{(n-3)!} \\
&=
\frac{n}{((n-1)!)^2(n-2)} ((n-1)t-(n-2)s)^{n-1}.
\end{align*}
Using the latter together with~\eqref{eq:proof-first-pair-inclusion-exclusion} and~\eqref{eq:proof-first-coordinate-volume}, we obtain
\begin{align*}
n!(n-1)!\operatorname{Vol}\bigl(\mathcal P(se_1+te_2)\bigr)
&=
n!(n-1)!
\left(
\frac1n V_1(s)+\frac1n\operatorname{Vol}\bigl(\mathcal R(s,t)\bigr)
\right)\\
&=
n^{n-1}s^{n-1}
+
\frac{n}{n-2}
\bigl((n-1)t-(n-2)s\bigr)^{n-1}.
\end{align*}
Theorem~\ref{thm:maj-conv-distr} now yields~\eqref{eq:maj-limit-two-point} whenever $\frac{n-2}{n-1}s\leq t\leq2s$.

It remains to extend the formula to arbitrary $s,t\geq0$.
Let $(\tilde{s},\tilde{t})$ be defined as in~\eqref{eq:first-adjacent-concavified-thresholds}.
By Lemma~\ref{lem:concave-ratios}, every $x\in\mathcal C_{n-1}$ satisfies
\begin{equation}
\frac{n-2}{n-1}x_1\leq x_2\leq2x_1.
\end{equation}
It follows that
\begin{equation}
\left\{
x\in\mathcal C_{n-1}
\colon
x_1\geq s,\ x_2\geq t
\right\}
=
\left\{
x\in\mathcal C_{n-1}
\colon
x_1\geq\tilde{s},\ x_2\geq\tilde{t}
\right\}.
\end{equation}
Since $(\zeta_1,\ldots,\zeta_{n-1})\in\mathcal C_{n-1}$ almost surely, we conclude that
\begin{equation}
\P(\zeta_1\geq s,\zeta_2\geq t)
=
\P(\zeta_1\geq\tilde{s},\zeta_2\geq\tilde{t}).
\end{equation}
The pair $(\tilde{s},\tilde{t})$ satisfies
\begin{equation}
\frac{n-2}{n-1}\tilde{s}
\leq
\tilde{t}
\leq
2\tilde{s},
\end{equation}
so the formula already proved under the assumption~\eqref{eq:proof-first-adjacent-compatible} applies to $(\tilde{s},\tilde{t})$.
\end{proof}

\begin{rem}
The one-point marginal in Corollary~\ref{cor:maj-limit-marginals} can also be obtained directly from the finite-sample marginal distribution.
Indeed, for every fixed $s\geq0$ and all sufficiently large $m$, formula~\eqref{eq:marginal_phik_maj} gives
\[
\P\left(
\phi_k(Y)
\leq
\frac{k}{n}+\frac{s}{m^{1/(n-1)}}
\right)
=
\frac{\mu_{n,k}s^{n-1}}{m}.
\]
It then follows from~\eqref{eq:distr_meet_maj} that
\[
\P\left(
m^{1/(n-1)}
\left(
\phi_k(\mathsf{\Lambda}_{n,m})-\frac{k}{n}
\right)
\geq s
\right)
=
\left(
1-\frac{\mu_{n,k}s^{n-1}}{m}
\right)^m
\xrightarrow{m\to\infty}
\exp\left(-\mu_{n,k}s^{n-1}\right),
\]
as required.
\end{rem}
In the following proofs, we shall repeatedly use the Bonferroni inequalities.
Given a finite family
of events $E_1,\ldots,E_n$, set
\[
\Sigma_r
:=
\sum_{1\le i_1<\cdots<i_r\le n}
\mathbb P(E_{i_1}\cap\cdots\cap E_{i_r}),
\qquad r\ge1.
\]
Then, for every $k\ge0$,
\begin{equation}
\label{eq:Bonferroni}
\sum_{r=1}^{2k} (-1)^{r-1} \Sigma_r
\leq \P\!\left(\bigcup_{i=1}^n E_i\right)
\leq \sum_{r=1}^{2k+1} (-1)^{r-1} \Sigma_r.
\end{equation}
\begin{proof}[Proof of Theorem~\ref{thm:maj-join-conv-distr}]
We divide the proof into four steps.

\medskip
\noindent
\emph{Step 1: lower tails of the sums of the smallest coordinates.} 
Let $X=(X_1,\ldots,X_n)$ be uniformly distributed on $\Delta_{n-1}$, and set $Y:=X^\downarrow$, so that $Y$ is uniformly distributed on $\Delta_{n-1}^{\downarrow}$.
For every $k\in\llbracket n-1\rrbracket$, define $T_k:=1-\phi_k(Y)$ to be the sum of the $q:=n-k$ smallest coordinates of $Y$.
Let us use the shorthands
\[
X_I:= \sum_{j\in I}X_j, \qquad
A_I(t)
:=
\left\{
X_I\leq t
\right\}.
\]
for each $I\subseteq\llbracket n\rrbracket$.
Thus,
\begin{equation}
\label{eq:Tk-tail-union}
\{T_k\leq t\}
=
\bigcup_{\substack{I\subseteq\llbracket n\rrbracket\\ |I|=q}}
A_I(t).
\end{equation}

We will now use \eqref{eq:Tk-tail-union} and the first truncation in the Bonferroni inequatlities~\eqref{eq:Bonferroni} to estimate $\mathbb{P}(T_k\leq t)$:
\begin{equation}\label{eq:BonferroniFirst}
    \Sigma_{1}-\Sigma_{2} \leq \mathbb{P}(T_k\leq t) \leq \Sigma_1,
\end{equation}
where 
\begin{equation*}
    \Sigma_1=\sum_{|I|=q}\mathbb{P}(A_I(t)), \qquad \Sigma_2=\sum_{\substack{I,J\subseteq \llbracket n \rrbracket\\|I|=|J|=q, I\neq J}} \mathbb P\bigl(A_I(t)\cap A_J(t)\bigr).
\end{equation*}
By Proposition~\ref{prop:UT-uniform}, $X$ has the same distribution as the vector of uniform spacings $S_n=(S_{1,n},\dots,S_{n,n})$.
In particular, its coordinates are exchangeable, and hence, for every $I\subseteq \llbracket n\rrbracket$ of cardinality $q$,
\[
X_I
\stackrel{\operatorname d}{=}
\sum_{j=1}^{q}S_{j,n}
=U_{q:n-1}
\stackrel{\operatorname d}{=}\phi_q(X).
\]
The distribution function of $\phi_q(X)$ from Lemma~\ref{lem:marginal_phik_UT} gives
\begin{align*}
\P(A_I(t))
=\sum_{i=q}^{n-1}
\binom{n-1}{i}t^i(1-t)^{n-1-i}
=\binom{n-1}{q}t^q+O(t^{q+1}),
\qquad
t\downarrow0.
\end{align*}

We next show that intersections corresponding to distinct subsets contribute only to the error term.
For any distinct $I,J\subseteq \llbracket n\rrbracket$ of cardinality $q$, set $d:=|I\setminus J|\geq1$.
On $A_I(t)\cap A_{J}(t)$,
\[
X_{I\cup J}
\leq
X_I+X_J
\leq2t.
\]
Since $|I\cup J|=q+d$, if $q+d<n$, exchangeability gives
\[
\P\bigl(A_I(t)\cap A_{J}(t)\bigr)
\leq
\P(\phi_{q+d}(X)\leq2t)
=O(t^{q+d}).
\]
If $q+d=n$, then $\sum_{j\in I\cup J}X_j=1$, so the intersection is empty for $t<1/2$.
Thus, in all cases,
\[
\P\bigl(A_I(t)\cap A_{J}(t)\bigr)
=
O(t^{q+1}).
\]

Since $n$ is fixed, there are only finitely many relevant subsets.
Then, applying the Bonferroni bounds~\eqref{eq:BonferroniFirst} to~\eqref{eq:Tk-tail-union} yields
\begin{equation}
\label{eq:Tk-tail}
\P(T_k\leq t)
=\binom{n}{q}\binom{n-1}{q}t^q + O(t^{q+1})
=\nu_{n,k}t^{n-k} + O(t^{n-k+1}),
\qquad
t\downarrow0,
\end{equation}
where $\nu_{n,k}$ is defined in~\eqref{eq:join-limit-marginals}.

\medskip
\noindent
\emph{Step 2: simultaneous rare events at different levels.}
Fix $1\leq k<\ell\leq n-1$ and write $q:=n-k$ and $r:=n-\ell$, so that $q>r$.
Let $I,J\subseteq\llbracket n\rrbracket$ be such that $|I|=q$ and $|J|=r$, and set $d:=|I\setminus J|$.
Since $q>r$, we have $d\geq q-r=\ell-k$.
By~\eqref{eq:Tk-tail-union}, we have
\begin{equation}
\label{eq:simultaneous-rare-events}
\{T_k\leq t, T_{\ell}\leq u\}
= \bigcup_{\substack{I,J\subseteq\llbracket n\rrbracket\\ |I|=q, \, |J|=r}}
A_J(u) \cap A_I(t)
\subseteq \bigcup_{\substack{I,J\subseteq\llbracket n\rrbracket\\ |I|=q, \, |J|=r}}
A_J(u) \cap A_{I\setminus J}(t),
\end{equation}
where the inclusion is due to the fact that the coordinates of $X$ are nonnegative.

We next show that, for all sufficiently small $t,u>0$,
\begin{equation}
\label{eq:simultaneous-rare-events-estimate}
\P(A_J(u) \cap A_{I\setminus J}(t))
\leq
C_{n,r,d}u^rt^d,
\end{equation}
where $C_{n,r,d}$ is a constant that depends on $n,r,d$ only.
If $r+d=n$, then $X_J+X_{I\setminus J}=1$, and the event $A_J(u) \cap A_{I\setminus J}(t)$ is empty for $t+u<1$.
From now on, we may thus suppose that $r+d\leq n-1$.
By the uniform-spacing representation and the exchangeability of the spacings,
\[
(X_J,X_{I\setminus J})
\stackrel{\operatorname d}{=}
\left(
U_{r:n-1},
U_{r+d:n-1}-U_{r:n-1}
\right).
\]
On the event
\[
\{U_{r:n-1}\leq u,
U_{r+d:n-1}-U_{r:n-1}\leq t\},
\]
there exist distinct labels $j_1,\ldots,j_r,i_1,\ldots,i_d$ such that $U_{j_r}\leq u$ and
\[
U_{j_a}\leq U_{j_r} \leq U_{i_b}\leq U_{j_r}+t \qquad \text{for all $a\in \llbracket r-1 \rrbracket$, $b\in \llbracket d \rrbracket$}.
\]
For each fixed choice of labels, condition on $U_{j_r}=x$.
The remaining variables being independent, the conditional probability that
\[
U_{j_a}\leq x \leq U_{i_b}\leq x+t \qquad \text{for all $a\in \llbracket r-1 \rrbracket$, $b\in \llbracket d \rrbracket$}
\]
is at most $x^{r-1}t^d$.
Integrating over $x\in[0,u]$ gives
\[
\int_0^u x^{r-1}t^d\,\diff x
=
\frac{u^rt^d}{r}.
\]
Since there are only finitely many choices of labels, a union bound proves~\eqref{eq:simultaneous-rare-events-estimate}.

Since $d\geq\ell-k$ and $t\leq1$, \eqref{eq:simultaneous-rare-events} and~\eqref{eq:simultaneous-rare-events-estimate} yield
\begin{equation}
\label{eq:two-level-negligible}
\P(T_k\leq t,\ T_\ell\leq u)
=
O\left(t^{\ell-k}u^{n-\ell}\right),
\qquad
t,u\downarrow0.
\end{equation}

\medskip
\noindent
\emph{Step 3: joint limit for the componentwise maxima.}
Let
\[
M_{m,k}
:=
\max_{i\in \llbracket m\rrbracket}\phi_k(Y_{n,i})
\]
and fix $s_1,\ldots,s_{n-1}\geq0$.
Since $Y_{n,1},Y_{n,2},\dots$ are i.i.d.\ copies of $Y$,
\[
\begin{aligned}
&\P\left(\bigcap_{k=1}^{n-1}
\left\{m^{1/(n-k)}(1-M_{m,k})>s_k\right\}\right)
=\P\left(\bigcap_{k=1}^{n-1}
\left\{\min_{i\in \llbracket m\rrbracket}\{1-\phi_k({Y_{n,i}})\}>\frac{s_k}{m^{1/(n-k)}}\right\}\right) \\
&=\P\left(\bigcap_{k=1}^{n-1}\left\{1-\phi_k({Y})>\frac{s_k}{m^{1/(n-k)}}\right\}\right)^m
=\left[
1-\P\left(
\bigcup_{k=1}^{n-1}E_{m,k}
\right)
\right]^m,
\end{aligned}
\]
where
\[
E_{m,k}
:=
\left\{
T_k\leq s_km^{-1/(n-k)}
\right\}.
\]

Taking $t=s_km^{-1/(n-k)}$ in~\eqref{eq:Tk-tail} gives
\[
\P(E_{m,k})
=
\nu_{n,k}\frac{s_k^{n-k}}{m} + o\left(\frac{1}{m}\right).
\]
Moreover, for $1\leq k<\ell \leq n-1$, taking $t=s_km^{-1/(n-k)}$ and $u=s_\ell m^{-1/(n-\ell)}$ in~\eqref{eq:two-level-negligible} gives
\begin{equation*}
\P\left(
E_{m,k} \cap E_{m,\ell}
\right)
=
O\left(
m^{-1-(\ell-k)/(n-k)}
\right)
= o\left(\frac{1}{m}\right).
\end{equation*}
The Bonferroni inequalities therefore yield
\[
\P\left(
\bigcup_{k=1}^{n-1}E_{m,k}
\right)
=
\sum_{k=1}^{n-1} \nu_{n,k}\frac{s_k^{n-k}}{m} + o\left(\frac{1}{m}\right).
\]

We conclude that
\[
\P\left(\bigcap_{k=1}^{n-1}
\left\{m^{1/(n-k)}(1-M_{m,k})>s_k\right\}\right)
\to
\exp\left(
-\sum_{k=1}^{n-1}\nu_{n,k}s_k^{n-k}
\right).
\]
The right-hand side is the joint survival function of independent random variables $\eta_1,\ldots,\eta_{n-1}$ such that
\[
\eta_k^{n-k}\sim\Exp(\nu_{n,k}).
\]
It follows that
\begin{equation}
\label{eq:Rmk-joint-limit}
\left(
m^{1/(n-k)}(1-M_{m,k})
\right)_{k=1}^{n-1}
\xrightarrow{\operatorname d}
(\eta_1,\ldots,\eta_{n-1}).
\end{equation}

\medskip
\noindent
\emph{Step 4: the componentwise maxima are asymptotically concave.}
We now prove~\eqref{eq:concavification-inactive}, i.e.,
\begin{equation}
\label{eq:concavification-inactive-proof}
\P\left(
\phi_k(\mathsf{V}_{n,m})=M_{m,k}
\text{ for all }k\in\llbracket n-1\rrbracket
\right)
\longrightarrow1
\end{equation}
as $m\to\infty$.
This, combined with~\eqref{eq:Rmk-joint-limit}, implies also~\eqref{eq:join-joint-limit}.

Set $M_{m,0}:=0$, $M_{m,n}:=1$, and
\[
R_{m,k}:=1-M_{m,k}, \qquad 0\leq k\leq n.
\]
The augmented sequence $(M_{m,0},M_{m,1},\ldots,M_{m,n})$ is concave if and only if
\begin{equation}
\label{eq:R-concave}
2R_{m,k}
\leq
R_{m,k-1}+R_{m,k+1},
\qquad
k\in\llbracket n-1\rrbracket,
\end{equation}
i.e., if and only if $(R_{m,0},R_{m,1},\ldots,R_{m,n})$ is convex.

Fix $k\in\llbracket2,n-1\rrbracket$.
By \eqref{eq:Rmk-joint-limit}, the continuous mapping theorem, and the fact that $\eta_{k-1}>0$ almost surely,
\[
\frac{
m^{1/(n-k)}R_{m,k}
}{
m^{1/(n-k+1)}R_{m,k-1}
}
\xlongrightarrow{\operatorname d}
\frac{\eta_k}{\eta_{k-1}}.
\]
By Slutsky's theorem,
\[
\frac{R_{m,k}}{R_{m,k-1}}
=
m^{-1/((n-k)(n-k+1))}
\frac{
m^{1/(n-k)}R_{m,k}
}{
m^{1/(n-k+1)}R_{m,k-1}
}
\xlongrightarrow{\mathrm{d}}0, \qquad k\in\llbracket2,n-1\rrbracket.
\]
The same holds for $k=1$: since $R_{m,0}=1$ and $R_{m,1}\to0$ in distribution by~\eqref{eq:Rmk-joint-limit}, we have $R_{m,1}/R_{m,0}\to 0$ in distribution.
Since $n$ is fixed, a union bound gives
\[
\P\left(
2R_{m,k}\leq R_{m,k-1}
\text{ for every }k\in\llbracket n-1\rrbracket
\right)
\longrightarrow1.
\]
Using $R_{m,k+1}\geq0$, we deduce that~\eqref{eq:R-concave} holds with probability tending to $1$ as $m\to\infty$.
Consequently,
\begin{equation}
\label{eq:Mmk-concavity}
\P\left(
(M_{m,0},M_{m,1},\ldots,M_{m,n})
\text{ is concave}
\right)
\longrightarrow1.
\end{equation}

By Proposition~\ref{prop:meet-majorisation},
\[
\widehat{\phi}(\mathsf{V}_{n,m})
=
(0,M_{m,1},\ldots,M_{m,n-1},1)^*.
\]
On the event in~\eqref{eq:Mmk-concavity}, the concave envelope leaves the augmented sequence unchanged, thus yielding~\eqref{eq:concavification-inactive-proof}.
\end{proof}

\subsubsection{High-dimensional results}

\begin{proof}[Proof of Theorem~\ref{thm:maj-meet-m-n-collapse}]
By~\eqref{eq:meet_maj_ineq},
\begin{equation}
\label{eq:meet-maj-infty-norm}
\norm{\mathsf{\Lambda}_{n,m_n}-\bot_n}_{\infty}
=
\max\left\{
\phi_1(\mathsf{\Lambda}_{n,m_n})-\frac1n,\,
\phi_{n-1}(\mathsf{\Lambda}_{n,m_n})-\frac{n-1}{n}
\right\}.
\end{equation}

For every $\epsilon\in(0,1)$, \eqref{eq:distr_meet_maj} and~\eqref{eq:marginal_n-1_lower_maj} give
\[
\P\left(
n\norm{\mathsf{\Lambda}_{n,m_n}-\bot_n}_{\infty}\geq\epsilon
\right)
\geq \P\left(\phi_{n-1}(\mathsf{\Lambda}_{n,m_n})\geq \frac{n-1}{n}+\frac{\epsilon}{n}\right)
=\left(1-\epsilon^{n-1}\right)^{m_n}.
\]
Suppose that $n\norm{\mathsf{\Lambda}_{n,m_n}-\bot_n}_{\infty}\to0$ almost surely, and hence in probability.
Then, $\left(1-\epsilon^{n-1}\right)^{m_n}\to 0$.
Taking logarithms, we obtain $m_n \log(1-\epsilon^{n-1}) \to-\infty$.
Since $\log(1-\epsilon^{n-1})\sim -\epsilon^{n-1}$ as $n\to\infty$, we have $m_n\epsilon^{n-1}\to\infty$.
Taking again logarithms, we see that $\log m_n +(n-1)\log\epsilon \to\infty$.
In particular, for all $\epsilon\in (0,1)$,
\begin{equation}
\label{eq:log-mn-infinity-proof}
\frac{\log m_n}{n} + \frac{n-1}{n}\log\epsilon>0
\qquad
\text{for all sufficiently large $n$}.
\end{equation}
If $\log m_n/n$ does not go to $\infty$, then there exist a positive constant $C$ such that $\frac{\log m_n}{n}\leq C$ along a subsequence.
Taking $\epsilon<e^{-C}$, we see that the sequence in~\eqref{eq:log-mn-infinity-proof} is negative for all $n$ along the subsequence, which is a contradiction.
Therefore, $\log m_n/n\to\infty$.

Conversely, suppose that $\log m_n/n\to\infty$.
By \eqref{eq:meet-maj-infty-norm},
\[
\P\left(
n\norm{\mathsf{\Lambda}_{n,m_n}-\bot_n}_{\infty}\geq\epsilon
\right)
\leq \P\left(\phi_1(\mathsf{\Lambda}_{n,m_n})\geq \frac{1+\epsilon}{n}\right)
+ \P\left(\phi_{n-1}(\mathsf{\Lambda}_{n,m_n})\geq \frac{n-1}{n}+\frac{\epsilon}{n}\right).
\]
We have an exact expression~\eqref{eq:marginal_n-1_lower_maj} for the second probability on the right-hand side, but not for the first one, since~\eqref{eq:marginal_1_lower_maj} does not apply.
Let us then find an alternative estimate for the first probability.
If $\epsilon\in (0,1)$, the intervals
\[
I_j
:=
\left[
\frac{j-\epsilon/2}{n},
\frac{j+\epsilon/2}{n}
\right],
\qquad
j\in\llbracket n-1\rrbracket
\]
are pairwise disjoint and contained in $[0,1]$.
Let $U_1,\ldots,U_{n-1}$ be independent uniform random variables on $[0,1]$.
If each interval $I_j$ contains exactly one of the $U_i$, then every uniform spacing is at most $(1+\epsilon)/n$.

Since the intervals are pairwise disjoint and have length $\epsilon/n$, the probability that each of them contains exactly one of the $U_i$ is $(n-1)!(\epsilon/n)^{n-1}$.
By Proposition~\ref{prop:maj-uniform}, if $Y$ is uniformly distributed on $\Delta_{n-1}^{\downarrow}$, then $\phi_1(Y)$ has the same distribution as the largest uniform spacing.
It follows that
\[
\P\left(
\phi_1(Y)\leq\frac{1+\epsilon}{n}
\right)
\geq
(n-1)!\left(\frac{\epsilon}{n}\right)^{n-1}
\geq
\left(\frac{\epsilon}{2e}\right)^{n-1},
\]
where we used $k!\geq(k/e)^k$ and $(n-1)/n\geq1/2$.
Consequently,
\[
\P\left(
\phi_1(\mathsf{\Lambda}_{n,m_n})\geq\frac{1+\epsilon}{n}
\right)
\leq
\left(
1-\left(\frac{\epsilon}{2e}\right)^{n-1}
\right)^{m_n}.
\]

Using $1-x\leq e^{-x}$, we obtain
\begin{align*}
\P\left(
n\norm{\mathsf{\Lambda}_{n,m_n}-\bot_n}_{\infty}\geq\epsilon
\right)
&\leq
\left(
1-\left(\frac{\epsilon}{2e}\right)^{n-1}
\right)^{m_n}
+
\left(1-\epsilon^{n-1}\right)^{m_n}
\\
&\leq
2\exp\left(
-m_n\left(\frac{\epsilon}{2e}\right)^{n-1}
\right).
\end{align*}
The right-hand side is summable in $n$; indeed,
\[
\frac1n
\log\left(
m_n\left(\frac{\epsilon}{2e}\right)^{n-1}
\right)
=
\frac{\log m_n}{n}
+
\frac{n-1}{n}\log\left(\frac{\epsilon}{2e}\right)
\longrightarrow\infty,
\]
and hence $m_n\left(\epsilon/(2e)\right)^{n-1}\geq e^n$ for large enough $n$.
The Borel--Cantelli lemma yields $n\norm{\mathsf{\Lambda}_{n,m_n}-\bot_n}_{\infty}
\xrightarrow{\operatorname{a.s.}}0$, as required.
\end{proof}

\begin{proof}[Proof of Theorem~\ref{thm:maj-join-m-n-scaling}]
Let
\[
S_n^{(i)}
=
\bigl(S_{1,n}^{(i)},\ldots,S_{n,n}^{(i)}\bigr),
\qquad i\geq1,
\]
be independent copies of the uniform-spacing vector $S_n=(S_{1,n},\dots,S_{n,n})$ (see Sec.~\ref{sec:spacings}).
By Proposition~\ref{prop:maj-uniform}, we may take
$Y_{n,i}=(S_n^{(i)})^\downarrow$.

We first observe that concavification does not change the first partial sum.
Indeed, if $M:=\max_{i\in\llbracket m_n\rrbracket}\phi_1(Y_{n,i})$, then $\phi_k(Y_{n,i})\leq kM$ for every $i$ and $k$.
The concave sequence $\bigl(\min\{kM,1\}\bigr)_{k=0}^n$ therefore equals $1$ at $k=n$ (as $M\geq 1/n$), dominates the componentwise maxima of the partial sums, and agrees with them at $k=1$.
The description of the majorisation join in Proposition~\ref{prop:meet-majorisation} consequently gives
\[
\phi_1(\mathsf{V}_{n,m_n})
=
\max_{i\in\llbracket m_n\rrbracket}\phi_1(Y_{n,i})
=
\max_{\substack{i\in\llbracket m_n\rrbracket\\j\in\llbracket n\rrbracket}}
S_{j,n}^{(i)}.
\]
Thus,
\[
D_n^\vee
=
1-
\max_{\substack{i\in\llbracket m_n\rrbracket\\j\in\llbracket n\rrbracket}}
S_{j,n}^{(i)}.
\]

Fix $\epsilon\in(0,1)$.
By Lemma~\ref{lem:joint-upper-tails-spacings}, each uniform spacing $S_{j,n}^{(i)}$ satisfies $\P(S_{j,n}^{(i)}>1-t)=t^{n-1}$ for $t\in[0,1]$.
Therefore, the union bound gives
\[
\P\left(D_n^\vee<(1-\epsilon)a_n\right)
\leq
nm_n\bigl((1-\epsilon)a_n\bigr)^{n-1}
=
(1-\epsilon)^{n-1},
\]
where $a_n:=(nm_n)^{-1/(n-1)}$ is defined in~\eqref{eq:maj-join-an}.
On the other hand, if $(1+\epsilon)a_n\leq1$, independence between the $m_n$ spacing vectors and the inequality $1-x\leq e^{-x}$ give
\[
\begin{aligned}
\P\left(D_n^\vee>(1+\epsilon)a_n\right)
&\leq
\P\left(
S_{1,n}^{(i)}
<
1-(1+\epsilon)a_n
\text{ for every }i\in\llbracket m_n\rrbracket
\right) \\
&=
\left(
1-\bigl((1+\epsilon)a_n\bigr)^{n-1}
\right)^{m_n}
\leq
\exp\left(-\frac{(1+\epsilon)^{n-1}}{n}\right).
\end{aligned}
\]
If $(1+\epsilon)a_n>1$, the bound still holds, since the probability on the left is zero.
Both bounds are summable in $n$, so the Borel--Cantelli lemma proves~\eqref{eq:maj-join-m-n-scaling}.

We now prove the Gumbel limit.
For $x\in\R$ and for all sufficiently large $n$,
\[
\P\left(
n\left(1-\frac{D_n^\vee}{a_n}\right)\leq x
\right)
=
\P\left(
\max_{j\in\llbracket n\rrbracket}S_{j,n}\leq 1-a_n\left(1-\frac{x}{n}\right)
\right)^{m_n}
= (1-p_n)^{m_n},
\]
where
\[
p_n:=\P\left(\bigcup_{j=1}^n E_j\right),
\qquad
E_j:= \left\{S_{j,n}>1-a_n\left(1-\frac{x}{n}\right)\right\}.
\]
We will repeatedly use the Bonferroni inequalities~\eqref{eq:Bonferroni}.

Suppose first that $m_n\to\infty$.
By Lemma~\ref{lem:joint-upper-tails-spacings}, for large enough $n$ we have
\[
\Sigma_1 = \sum_{j=1}^n \P(E_j)
= n\left(a_n\left(1-\frac{x}{n}\right)\right)^{n-1}
= \frac{e^{-x}}{m_n} +o\left(\frac{1}{m_n}\right)
\]
and, using the inequality $2y-1\leq y^2$ for $y\in\R$,
\[
\begin{aligned}
\Sigma_2
&=\sum_{1\leq i<j\leq n} \P(E_i\cap E_j)
\leq \binom{n}{2} \left(2a_n\left(1-\frac{x}{n}\right)-1\right)^{n-1} \\
&\leq n^2 \left(a_n\left(1-\frac{x}{n}\right)\right)^{2(n-1)}
\leq \frac{1}{m_n^2} \left(1-\frac{x}{n}\right)^{2(n-1)}
= o\left(\frac{1}{m_n}\right).
\end{aligned}
\]
It thus follows from the Bonferroni inequalities $\Sigma_1-\Sigma_2\leq p_n \leq \Sigma_1$ that
\[
(1-p_n)^{m_n} = \left(1- \frac{e^{-x}}{m_n} + o\left(\frac{1}{m_n}\right)\right)^{m_n} \xlongrightarrow{n\to\infty} \exp\left(-e^{-x}\right).
\]

Suppose next that $m_n= m$ for some fixed positive integer $m$, so that
\begin{align*}
a_n
&= \exp\left(-\frac{\log(nm)}{n-1}\right)
= 1- \frac{\log(nm)}{n-1} +O\left(\frac{(\log n)^2}{n^2}\right), \\
a_n\left(1-\frac{x}{n}\right)
&= 1-\frac{\log(nm)}{n-1}-\frac{x}{n}+ O\left(\frac{(\log n)^2}{n^2}\right).
\end{align*}
Again by Lemma~\ref{lem:joint-upper-tails-spacings}, we have, for every fixed $r\geq 1$ and all sufficiently large $n$,
\[
\begin{aligned}
\Sigma_r
&= \binom{n}{r}\left(1-r\left(1-a_n\left(1-\frac{x}{n}\right)\right)\right)^{n-1} \\
&\sim \frac{n^r}{r!}
\exp\left\{(n-1)\log\left(1-r\left(\frac{\log(nm)}{n-1}+\frac{x}{n}+ O\left(\frac{(\log n)^2}{n^2}\right)\right)\right)\right\} \\
&= \frac{n^r}{r!}
e^{-r(\log n + \log m + x) +o(1)}
\xlongrightarrow{n\to\infty}
\frac{e^{-rx}}{m^r r!}.
\end{aligned}
\]
Thus, letting first $n\to\infty$ and then $k\to\infty$ in~\eqref{eq:Bonferroni}, we obtain
\[
p_n \longrightarrow -\sum_{r=1}^{\infty} \frac{(-1)^{r}}{r!} \left(\frac{e^{-x}}{m}\right)^r
= 1-\exp\left(-\frac{e^{-x}}{m}\right).
\]
We conclude that $(1-p_n)^m\to \exp\left(-e^{-x}\right)$.

Finally, let $m_n$ be any integer sequence.
Every subsequence contains a further subsequence along which either $m_n\to\infty$ or $m_n$ is constant.
The preceding two cases therefore show that
\[
\P\left(
n\left(1-\frac{D_n^\vee}{a_n}\right)\leq x
\right)
\longrightarrow
\exp\left(-e^{-x}\right),
\qquad x\in\R.
\]
This is~\eqref{eq:maj-join-m-n-Gumbel}.
\end{proof}

\begin{proof}[Proof of Corollary~\ref{cor:maj-join-m-n-collapse}]
By~\eqref{eq:join_maj_ineq}, for every $p\in[1,\infty]$,
\[
D_n^\vee
=
\bigl\|\mathsf{V}_{n,m_n}-\top_n\bigr\|_\infty
\leq
\bigl\|\mathsf{V}_{n,m_n}-\top_n\bigr\|_p
\leq
\bigl\|\mathsf{V}_{n,m_n}-\top_n\bigr\|_1
=
2D_n^\vee.
\]
Moreover,
\[
\log a_n=-\frac{\log n+\log m_n}{n-1},
\]
so $\log m_n/n\to c\in[0,\infty]$ if and only if $a_n\to e^{-c}$.
The conclusions now follow from~\eqref{eq:maj-join-m-n-scaling} exactly as in the proof of Corollary~\ref{cor:UT-m-n-collapse}.
\end{proof}

\section{Future directions}
\label{sec:future}

The results of this paper suggest several directions for further investigation.

\subsection{Explicit structure of the majorisation-meet limit}

At fixed dimension, the joint distribution of the limiting majorisation meet is described in terms of the volumes $\operatorname{Vol}\left(\mathcal P(s_1,\ldots,s_{n-1})\right)$; 
 see~\eqref{eq:maj-volume-formula}. 
The inclusion--exclusion representation obtained in Proposition~\ref{prop:maj-volume-piecewise} shows that these volumes are piecewise polynomial functions of the truncation parameters.
Although this gives explicit one-point marginals and some low-dimensional marginals, the formulae become increasingly complicated as the number of coordinates grows.

It would be useful to describe the chamber decomposition of the parameter space on which these volumes are polynomial and to develop an efficient method for computing the corresponding polynomial pieces.
Such a description could yield more tractable formulae for arbitrary finite-dimensional marginals and clarify the dependence structure of the limiting vector.
A direct probabilistic representation or sampling construction for this vector would also be valuable and might provide a more conceptual explanation of its palindromic symmetry.

\subsection{High-dimensional profiles and the exponential regime}

For the majorisation meet, a particularly natural unresolved regime is
\[
\frac{\log m_n}{n}\longrightarrow c\in(0,\infty).
\]
The collapse theorem identifies the superexponential regime \(c=\infty\), but does not describe the limiting behaviour of
\[
n\bigl\|\mathsf{\Lambda}_{n,m_n}-\bot_n\bigr\|_\infty
\]
when \(c\) is finite.
By~\eqref{eq:meet_maj_ineq}, this norm is determined by the first and last nontrivial partial sums.
The latter has an explicit distribution, whereas the former satisfies
\[
\phi_1(\mathsf{\Lambda}_{n,m_n})
=
\min_{i\in\llbracket m_n\rrbracket}
\max_{j\in\llbracket n\rrbracket}S_{j,n}^{(i)}.
\]
Consequently, for \(a>0\),
\[
\P\left(
\phi_1(\mathsf{\Lambda}_{n,m_n})>\frac{1+a}{n}
\right)
=
\left[
1-
\P\left(
\max_{j\in\llbracket n\rrbracket}S_{j,n}
\leq\frac{1+a}{n}
\right)
\right]^{m_n}.
\]
Since \(S_n\) is uniformly distributed on \(\Delta_{n-1}\), the probability inside the brackets is
\[
\frac{
\operatorname{Vol}\left(
\left\{
x\in\Delta_{n-1}\colon
\max_{j\in\llbracket n\rrbracket}x_j\leq\frac{1+a}{n}
\right\}
\right)
}{
\operatorname{Vol}(\Delta_{n-1})
}.
\]
Sharp asymptotics for this volume should therefore lead to a description of the limiting behaviour and fluctuations of
\(n\|\mathsf{\Lambda}_{n,m_n}-\bot_n\|_\infty\) in the exponential regime.

More generally, the high-dimensional results obtained here principally describe distances from the extremal elements and, for the majorisation join, the first partial sum.
A natural next step is to study the full partial-sum profiles, for example when \(k/n\to\alpha\in(0,1)\).
For the join, this would require understanding the concavification of the componentwise maxima simultaneously over all coordinates and determining whether it is asymptotically inactive, or otherwise describing its effect on the limiting profile.

\subsection{Other distributions and other orders}

In this paper, we have focused on points sampled uniformly from \(\Delta_{n-1}\) and \(\Delta_{n-1}^{\downarrow}\).
The same questions can be investigated for more general Dirichlet distributions, including families whose parameters vary with the dimension, as well as for random-matrix-induced measures arising in quantum entanglement theory~\cite{Cunden20,Cunden2021,jain2024entangled}.
The behaviour of the density near the boundary of the simplex may change the relevant tail exponents, scaling regimes, and collapse thresholds.
It would be interesting to determine which of the phenomena found here are universal and which depend on the particular sampling measure.

One can also study joins and meets of independent random elements in other lattices.
Natural examples include probability simplices equipped with more general semigroup majorisation relations~\cite{MOA79,Cunden25,Cunden26}, as well as thermomajorisation orders~\cite{Chitambar19,vomEnde}.
The geometry of the underlying order should again determine the scaling of the extremal fluctuations, but new features may arise when the lattice operations are not described by ordinary coordinatewise extrema or concavification.

\section*{Statements and acknowledgements}

\subsection*{Competing interests}
\addtocontents{toc}{\protect\setcounter{tocdepth}{1}}
The authors have no relevant financial or non-financial interests to disclose.

\subsection*{Data Availability}
All numerical results presented in this paper were generated through simulations.
The corresponding simulation code is available from the authors upon reasonable request.

\subsection*{AI use disclosure}
The authors used ChatGPT (OpenAI; mainly model GPT-5.6 Sol) for mathematical and editorial assistance during the preparation of this manuscript.
The authors conjectured the fixed-dimensional scaling limit for the majorisation join stated in Theorem~\ref{thm:maj-join-conv-distr}, and ChatGPT generated a complete proof.
ChatGPT also proposed the inclusion--exclusion formula in Proposition~\ref{prop:maj-volume-piecewise} and the Gumbel limits in Theorems~\ref{thm:UT-m-n-scaling} and~\ref{thm:maj-join-m-n-scaling} and generated complete proofs.
The authors independently verified these arguments, revised their presentation, and take full responsibility for all results and proofs in the manuscript.
ChatGPT was also used to simplify proofs, and assist with editing and proofreading.
All images were created using Mathematica and later processed with Claude Sonnet 5.5 for visual rendering.

\subsection*{Funding}
This work is supported by the Italian National Group for Mathematical Analysis, Probability and their Applications (INdAM-GNAMPA), the Italian National Group of Mathematical Physics (INdAM-GNFM), the Istituto Nazionale di Fisica Nucleare (INFN) through the project ``QUANTUM'', the PRIN 2022 project 2022TEB52W-PE1 - ``The charm of integrability: from nonlinear waves to random matrices'',  the PNRR MUR project CN00000013 ``Italian National Centre on HPC, Big Data and Quantum Computing'',  the PNRR MUR project PE0000023-NQSTI, and the project ``Patto territoriale sistema universitario pugliese''.

\addtocontents{toc}{\protect\setcounter{tocdepth}{2}}

\appendix

\section{Geometry of the UT-majorisation and majorisation lattices}
\label{sec:geometry}

\subsection{Proofs of the lattice constructions}
\label{subsec:geometry-partial-sums}

\begin{proof}[Proof of Proposition~\ref{prop:meet-UT}]
We prove the identity for the meet, as the proof for the join is symmetric.
Set $u:=\phi(x)$, $v:=\phi(y)$, and $r:=\min(u,v)$.
Since $r$ is nondecreasing and has last coordinate equal to $1$, Lemma~\ref{lem:bijection} yields a unique $z\in\Delta_{n-1}$ such that $\phi(z)=r$.

Since $r\leq u$ and $r\leq v$, the element $z$ is a common lower bound of $x$ and $y$.
Conversely, if $w$ is any common lower bound, then $\phi(w)\leq u$ and $\phi(w)\leq v$, and hence $\phi(w)\leq r=\phi(z)$.
Thus $w\preceq z$, proving that $z=x\wedge y$.
\end{proof}

\begin{proof}[Proof of Proposition~\ref{prop:meet-majorisation}]
Set $u:=\widehat{\phi}(x)$ and $v:=\widehat{\phi}(y)$.
Both sequences are nondecreasing and concave, and satisfy $u_0=v_0=0$ and $u_n=v_n=1$.
Their componentwise minimum $r:=\min(u,v)$ is nondecreasing, has endpoints $0$ and $1$, and is also concave.
Indeed, if $r_k=u_k$, then
\[
r_{k-1}+r_{k+1}
\leq u_{k-1}+u_{k+1}
\leq 2u_k
=2r_k,
\]
and a similar argument applies if $r_k=v_k$.
By Lemma~\ref{lem:bijection}, $r$ is the (augmented) partial-sum sequence of an element of $\Delta_{n-1}^{\downarrow}$.
The argument in Proposition~\ref{prop:meet-UT} then shows that this element is $x\wedge y$.

The componentwise maximum $s:=\max(u,v)$ need not be concave (see Example~\ref{ex:concave} below), so we replace it by its concave envelope $s^*$.
This is again the partial-sum sequence of an element of $\Delta_{n-1}^{\downarrow}$.
Indeed, the concave sequence $(0,1,\ldots,1)$ is clearly a componentwise majorant of $s$.
Therefore, by the minimality of the concave envelope,
\[
s\leq s^*\leq(0,1,\ldots,1).
\]
Thus $s^*$ has endpoints $0$ and $1$; moreover, concavity forces it to be nondecreasing, since a negative increment would make all subsequent increments negative and contradict $s_n^*=1$.

Since $s^*\geq u$ and $s^*\geq v$, the corresponding element in $\Delta_{n-1}^{\downarrow}$ is a common upper bound of $x$ and $y$.
Conversely, the partial-sum sequence of any common upper bound is a concave majorant of $s$ and therefore, by the minimality of the concave envelope, lies above $s^*$.
Hence the corresponding element is $x\vee y$.
\end{proof}

\begin{proof}[Proof of Theorem~\ref{thm:lattices}]
For every $x\in\Delta_{n-1}^{\downarrow}$ and $k\in\llbracket n\rrbracket$, the average of the $k$ largest components of $x$ is at least the average of all its components.
Hence
\[
\phi_k(x)\geq\frac{k}{n}
=\phi_k\left(\bot_n\right),
\]
i.e., $\bot_n\preceq x$ for all $x\in\Delta_{n-1}^{\downarrow}$.
Moreover, for all $x\in\Delta_{n-1}$ (in particular, for $x\in\Delta_{n-1}^{\downarrow}$), $\phi_k(x)\leq1=\phi_k(\topUT_n)$ ($k\in \llbracket n-1\rrbracket$), so $x\preceq\topUT_n=\top_n$.
Finally, for all $x\in \Delta_{n-1}$, $\phi_k(x)\geq 0=\phi_k(\botUT_n)$ ($k\in \llbracket n-1\rrbracket$), so $\botUT_n\preceq x$.
This proves~\eqref{eq:UT_top_bottom} and~\eqref{eq:maj_top_bottom}.
The meet and join constructions of Propositions~\ref{prop:meet-UT} and~\ref{prop:meet-majorisation} complete the proof that both posets are bounded lattices.
\end{proof}

\subsection{Lorenz curves and matrix characterisations}

The UT-majorisation and majorisation orders admit geometric interpretations in terms of Lorenz curves, as well as characterisations in terms of stochastic matrix transformations.
The Lorenz-curve picture also illustrates the constructions of their meets and joins.

\begin{definition}
The \textbf{Lorenz curve} $L_x$ associated to $x\in\Delta_{n-1}$ is the piecewise affine function $L_x\colon [0,1]\to[0,1]$ interpolating the points
\[
(0,0),\left(\frac{1}{n},\phi_1(x)\right),\left(\frac{2}{n},\phi_2(x)\right),\ldots,\left(\frac{n-1}{n},\phi_{n-1}(x)\right),(1,1).
\]
\end{definition}

The function $L_x$ is nondecreasing, and the invertibility of $\phi$ makes the correspondence $x\leftrightarrow L_x$ one-to-one.
Moreover, $x\in\Delta_{n-1}^{\downarrow}$ if and only if $L_x$ is concave.
If, instead, one considers increasingly ordered vectors, the corresponding Lorenz curves are convex, in accordance with the usual convention~\cite{lorenz1905methods}.

\begin{prop}
Let $x,y\in\Delta_{n-1}$.
The following are equivalent:
\begin{enumerate}
    \item\label{item:UTmaj} $x\preceq y$.
    \item\label{item:UT_Lorenz} $L_x(t)\leq L_y(t)$ for all $t\in[0,1]$.
    \item\label{item:upperTriang} $x=y P$, for some upper-triangular row stochastic matrix $P$.
\end{enumerate}
Moreover, for every $k\in\{0,\ldots,n\}$, we have
\begin{align}
\label{eq:Lorenz_meet_UT}
L_{x\wedge y}\left(\frac{k}{n}\right)
&=
\min\left(
L_x\left(\frac{k}{n}\right),
L_y\left(\frac{k}{n}\right)
\right),\\
\label{eq:Lorenz_join_UT}
L_{x\vee y}\left(\frac{k}{n}\right)
&=
\max\left(
L_x\left(\frac{k}{n}\right),
L_y\left(\frac{k}{n}\right)
\right)
\end{align}
and, between consecutive grid points, the Lorenz curves of the meet and join are obtained by affine interpolation.
\end{prop}
The equivalence \ref{item:UTmaj}$\iff$\ref{item:upperTriang} is discussed by Parker and Ram~\cite{ParkerRam} and explains why $\preceq$ is called both upper-triangular (shortened UT) majorisation and `unordered majorisation'~\cite[Chap.~14]{MOA79}.
The equivalence \ref{item:UTmaj}$\iff$\ref{item:UT_Lorenz} follows from Definition~\ref{def:UT-majorisation}, and provides a simple pictorial condition for this relation.
The characterisation of the Lorenz curves of the meet and join follows from Proposition~\ref{prop:meet-UT}.
Notice that the minimum and maximum in~\eqref{eq:Lorenz_meet_UT}--\eqref{eq:Lorenz_join_UT} are taken at the grid points: the affine interpolants need not coincide with the pointwise minimum and maximum between grid points.

\begin{prop}
Let $x,y\in\Delta_{n-1}^{\downarrow}$.
The following are equivalent:
\begin{enumerate}
    \item\label{item:maj} $x\preceq y$.
    \item\label{item:maj_Lorenz} $L_x(t)\leq L_y(t)$ for all $t\in[0,1]$.
    \item\label{item:doublyStoch} $x=y P$, for some doubly stochastic matrix $P$.
\end{enumerate}
Moreover, for every $k\in\{0,\ldots,n\}$,
\begin{equation}
\label{eq:Lorenz_meet_maj}
L_{x\wedge y}\left(\frac{k}{n}\right)=\min\left\{L_x\left(\frac{k}{n}\right),L_y\left(\frac{k}{n}\right)\right\},
\end{equation}
while, with $f^*$ denoting the least concave majorant of a function $f$ on $[0,1]$,
\begin{equation}
\label{eq:Lorenz_join_maj}
L_{x\vee y}=\bigl(\max(L_x,L_y)\bigr)^*.
\end{equation}
\end{prop}
The equivalence~\ref{item:maj}$\iff$\ref{item:maj_Lorenz} follows from Definition~\ref{def:Majorisation}.
The equivalence~\ref{item:maj}$\iff$\ref{item:doublyStoch} is a celebrated theorem of Hardy, Littlewood and P\'olya~\cite{HLP52}.
The characterisation of the Lorenz curves of the meet~\eqref{eq:Lorenz_meet_maj} and join~\eqref{eq:Lorenz_join_maj} follows from Proposition~\ref{prop:meet-majorisation}.

\begin{figure}
\centering
\includegraphics[width=.45\textwidth]{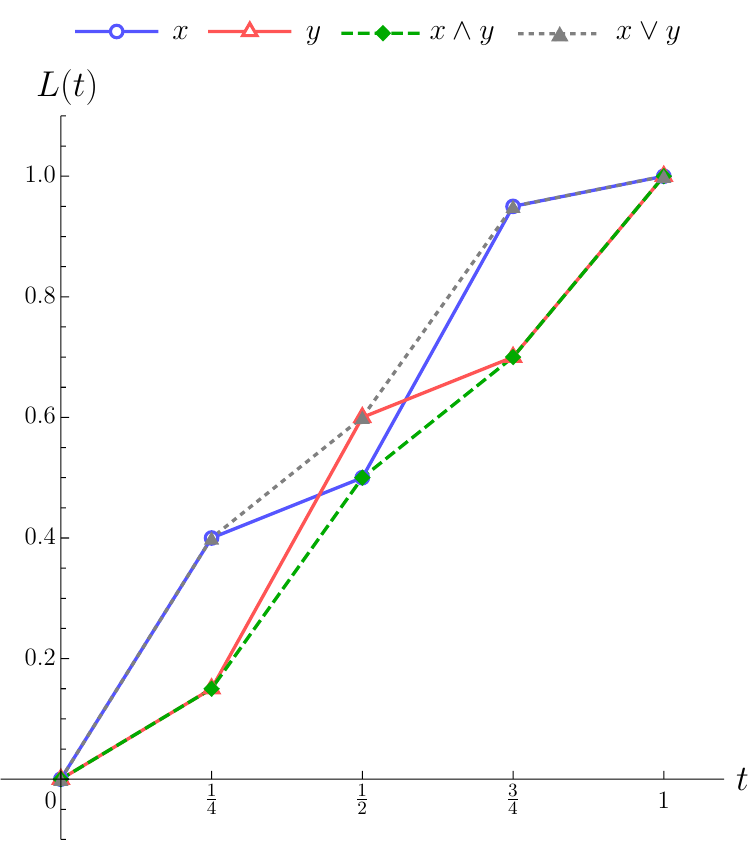}\quad 
\includegraphics[width=.45\textwidth]{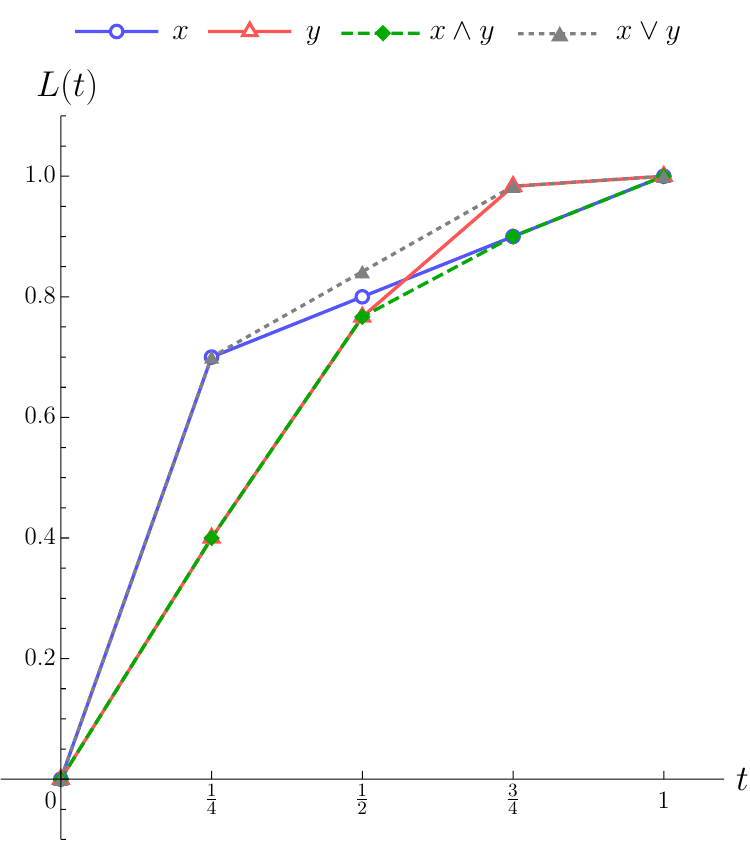}
        \caption{
        Lorenz curves of two incomparable vectors and their meet and join in the UT-majorisation lattice (left; Example~\ref{ex:ut-nonconcave-lorenz}) and the majorisation lattice (right; Example~\ref{ex:concave}). In the right-hand panel, the curve of the join lies above the gridwise maximum at the second interior grid point because concavification is required.}
        \label{fig:Lorenz}
    \end{figure}

\begin{ex}\label{ex:ut-nonconcave-lorenz}
For UT-majorisation, the Lorenz curves are nondecreasing but not necessarily concave.
Consider
\begin{equation}\label{eq:UTexample}
x=\left(
\frac{2}{5},
\frac{1}{10},
\frac{9}{20},
\frac{1}{20}
\right),
\qquad
y=\left(
\frac{3}{20},
\frac{9}{20},
\frac{1}{10},
\frac{3}{10}
\right).
\end{equation}
Their augmented partial-sum sequences are
\[
\widehat{\phi}(x)
=\left(
0,\frac{2}{5},\frac{1}{2},\frac{19}{20},1
\right), \qquad
\widehat{\phi}(y)
=\left(
0,\frac{3}{20},\frac{3}{5},\frac{7}{10},1
\right).
\]
In particular, $x$ and $y$ are incomparable, since
\[
\phi_1(x)>\phi_1(y),
\qquad
\phi_2(x)<\phi_2(y),
\qquad
\phi_3(x)>\phi_3(y).
\]
The UT-majorisation formulae give
\[
\widehat{\phi}(x\wedge y)
=
\left(
0,\frac{3}{20},\frac{1}{2},\frac{7}{10},1
\right), \qquad
\widehat{\phi}(x\vee y)
=
\left(
0,\frac{2}{5},\frac{3}{5},\frac{19}{20},1
\right).
\]
Equivalently,
\[x\wedge y
=
\left(
\frac{3}{20},
\frac{7}{20},
\frac{1}{5},
\frac{3}{10}
\right),
\qquad
x\vee y
=
\left(
\frac{2}{5},
\frac{1}{5},
\frac{7}{20},
\frac{1}{20}
\right).
\]
The gridwise minimum and maximum of $L_x$ and $L_y$ determine $L_{x\wedge y}$ and $L_{x\vee y}$.
See the left-hand panel of Figure~\ref{fig:Lorenz}.
\end{ex}

\begin{ex}\label{ex:concave}
For majorisation, the componentwise maximum of two admissible partial-sum sequences need not be concave.
For example, let
\begin{equation}
x=\left(\frac{7}{10},\frac{1}{10},\frac{1}{10},\frac{1}{10}\right),
\qquad
y=\left(\frac{2}{5},\frac{11}{30},\frac{13}{60},\frac{1}{60}\right)
\label{eq:xy_example}
\end{equation}
in $\Delta_3^\downarrow$.
Their augmented partial-sum sequences are
\[
\widehat{\phi}(x)
=
\left(0,\frac{7}{10},\frac{4}{5},\frac{9}{10},1\right), \qquad
\widehat{\phi}(y)
=
\left(0,\frac{2}{5},\frac{23}{30},\frac{59}{60},1\right).
\]
Their componentwise maximum is
\[
a
:=\max(\widehat{\phi}(x),\widehat{\phi}(y))
=\left(0,\frac{7}{10},\frac{4}{5},\frac{59}{60},1\right).
\]
This sequence is not concave, since
\[
2a_2=\frac{8}{5}
<
\frac{7}{10}+\frac{59}{60}
=\frac{101}{60}
=a_1+a_3.
\]
Its concave envelope is
\[
a^*
=
\left(0,\frac{7}{10},\frac{101}{120},\frac{59}{60},1\right)=\widehat{\phi}(x\vee y),
\]
so that
\[
x\vee y
=
\left(\frac{7}{10},\frac{17}{120},\frac{17}{120},\frac{1}{60}\right).
\]
In particular,
\[
\widehat{\phi}(x\vee y)
\neq
\max(\widehat{\phi}(x),\widehat{\phi}(y)).
\]
Thus the concavification in the construction of the majorisation join is genuinely necessary.
See the right-hand panel of Figure~\ref{fig:Lorenz}.
\end{ex}

\printbibliography

\end{document}